\documentclass[preprint,11pt]{elsarticle}
\usepackage{amsmath,amscd,amsthm,amsxtra,amsfonts,amssymb,algorithm,algorithmicx,algpseudocode,array,arydshln,arydshln,bbm,bm,booktabs,cases,caption,color,enumerate,extarrows,float,graphicx,geometry,hyperref,listings,mathrsfs,manyfoot,multirow,subfig,textcomp,verbatim,xcolor}

\usepackage{nomencl}
\usepackage{framed}
\usepackage{multicol}
\makenomenclature

\allowdisplaybreaks

\usepackage[title]{appendix}%

\numberwithin{equation}{section}

\newtheorem{theorem}{Theorem}
\newtheorem{lemma}[theorem]{Lemma}
\newtheorem*{lemma*}{Lemma}

\theoremstyle{definition}

\newtheorem{remark}[theorem]{Remark}

\begin{document}
	
	\begin{frontmatter}
		
		\title{The memory-dependent FPK equation for fractional Gaussian noise
		} 
		
\author[1,2]{Lifang Feng}\ead{FLF.fenglifang@outlook.com}

\author[1,2]{Bin Pei}\ead{binpei@nwpu.edu.cn}	

\author[1]{Yong Xu \corref{cor1}}\ead{hsux3@nwpu.edu.cn}

\affiliation[1]{organization={School of Mathematics and Statistics \& MOE Key Laboratory for Complexity Science in Aerospace, Northwestern Polytechnical University},
	city={Xi'an},
	postcode={710072},
	country={China}}
	
		\affiliation[2]{organization={	Research \& Development Institute of Northwestern Polytechnical University in Shenzhen}, 
		city={Shenzhen},
		postcode={518057},
		country={China}}

\cortext[cor1]{Corresponding author at: School of Mathematics and Statistics, Northwestern Polytechnical University, Xi'an, 710072, China.}

\begin{abstract}
This paper aims to explore non-Markovian dynamics of nonlinear dynamical systems subjected to fractional Gaussian noise (FGN) and Gaussian white noise (GWN). A novel memory-dependent Fokker-Planck-Kolmogorov (memFPK) equation is developed to characterize the probability structure in such non-Markovian systems. The main challenge in this research comes from the long-memory characteristics of FGN. These features make it impossible to model the FGN-excited nonlinear dynamical systems as finite dimensional GWN-driven Markovian augmented filtering systems, so the classical FPK equation is no longer applicable. To solve this problem, based on fractional Wick-It\^o-Skorohod integral theory, this study first derives the fractional It\^o formula. Then, a memory kernel function is constructed to reflect the long-memory characteristics from FGN. By using fractional It\^o formula and integration by parts, the memFPK equation is established. {Importantly, the proposed memFPK equation is not limited to specific forms of drift and diffusion terms, making it broadly applicable to a wide class of nonlinear dynamical systems subjected to FGN and GWN.} Due to the historical dependence of the memory kernel function, a Volterra adjustable decoupling approximation is used to reconstruct the memory kernel dependence term. This approximation method can effectively solve the memFPK equation, thereby obtaining probabilistic responses of nonlinear dynamical systems subjected to FGN and GWN excitations. Finally, some numerical examples verify the accuracy and effectiveness of the proposed method.

\end{abstract}
		
		
\end{frontmatter}

\begin{sloppypar}
\section{Introduction}\label{sec1}
Stochastic structural excitations are crucial in mechanical and vibration engineering \cite{li2009stochastic,li2012stochastic,cai2016elements,sun2006stochastic}. These excitations come from various sources, including machinery operational noise, unpredictable environmental forces like wind and earthquakes, and small manufacturing process variations \cite{liu2021generalized,obukhov2021modeling,vere1970stochastic}. Due to the stochastic and unpredictable nature of these excitations, it is essential to quantify structural response uncertainties from a probabilistic perspective. In many engineering cases, stochastic excitation goes beyond white noise \cite{molz1997fractional,einchcomb1995use,hu2016inverse}, making system responses depend on both current and historical states. Such memory effects make traditional Markovian models unsuitable for prediction. Here, non-Markovian modeling is necessary to accurately capture the system's memory characteristics. A typical example is fractional Gaussian noise (FGN) \cite{mandelbrot1968fractional,wang2022anomalous,jeon2010fractional,flandrin1989spectrum}, which is a non-white stochastic excitation with long memory and strong spatial-temporal correlations \cite{song2019fractional,li2010fractal,qian2003fractional}.
This study focuses on nonlinear dynamical system with long-memory properties and builds a non-Markovian system under the combined excitations of FGN and Gaussian white noise (GWN) to study non-Markovian dynamics. 
		
To trace the propagation of uncertainty and capture the probability information of stochastic dynamic responses, a relatively comprehensive theoretical framework has been established \cite{li2009stochastic,sun2006stochastic}, among which the Fokker-Planck-Kolmogorov (FPK) equation \cite{risken1996fokker} provides the most rigorous foundation. Under Markovian assumptions, the FPK framework has been widely used to analyze the probabilistic responses of stochastic dynamical systems, particularly for GWN, Poisson white noise \cite{chen2024probabilistic} and L\'evy noise \cite{sun2012fokker,sun2017governing,zan2020stochastic} excitations. 
{However, the classical FPK equation has limitations for non-Markovian stochastic systems with long-memory effects, such as those driven by FGN or involving fractional derivative operators. It is noteworthy that based on Kramers-Moyal expansion, the globally-evolving-based generalized density evolution equation (GE-GDEE)  \cite{li2006probability,lyu2022unified,lyu2024decoupled} method has significantly extended the applicability of probabilistic analysis to nonlinear dynamical systems endowed with fractional derivative elements subjected to GWN \cite{luo2022stochastic,luo2023equation}. Unfortunately, the intrinsic diffusion coefficients for FGN-driven systems cannot be directly estimated from the second-order conditional derivative moments within the GE-GDEE framework. Consequently, it is impossible to directly derive the system's response probability density functions (PDFs) using the established methodologies.} Therefore, developing probabilistic techniques to govern the evolution of the response PDFs for FGN-excited systems is critical for advancing the field of stochastic engineering dynamics.

So far, the research on probabilistic responses of nonlinear stochastic systems driven by FGN is still in its infancy, and the existing methods often have obvious limitations.  For instance, in relatively simple scenarios like linear systems under additive FGN excitation, Vyoral \cite{vyoral2005kolmogorov} established the corresponding Kolmogorov backward equation and then derived the associated FPK type equation. Additionally, Choi \cite{choi2021entropy} and Vas'ovskii \cite{vaskovskii2022analog} separately obtained the FPK type equation for nonlinear stochastic systems driven by FGN by using Doss-Sussman transformation technique and stochastic flow method, respectively. However, these studies failed to consider scenarios involving drift terms, thus severely limiting their applicability in engineering dynamics. For more general nonlinear systems-including those with external excitations, especially multiplicative excitation systems where randomness originates from the system's inherent structure-these challenges remain unresolved. In this regard, based on fractional Wick-It\^o-Skorohod (FWIS) integration and rough path theory, Pei et al. \cite{pei2025non}  have yielded memory-dependent probability density evolution equations (PDEEs)  for nonlinear systems excited by combined FGN and GWN under commutativity conditions, employing fractional power time-varying terms $ Ht^{2H-1} $ to capture the memory characteristics of stochastic excitations. However, the approach developed in \cite{pei2025non}  remain constrained by specific system coefficients with commutative conditions setting and cannot yet address general nonlinear systems. To overcome these limitations, this paper advances the memory-dependent PDEEs approach in \cite{pei2025non} by deriving a memory-dependent Fokker-Planck-Kolmogorov (memFPK) equation for general nonlinear dynamical system subjected to combined FGN and GWN without commutativity conditions. {A fundamental novelty of the proposed framework is the introduction of a memory-dependent kernel function (a Malliavin derivative term) via the FWIS integration theory. This kernel function systematically encodes the historical state information of the system. By utilizing the fractional It\^o formula and integration by parts, the resulting memFPK equation provides a probabilistic characterization of the system dynamics. Importantly, unlike many existing methods as mentioned above, the memFPK equation is not limited to specific forms of drift or diffusion terms, making it broadly applicable to a wide class of non-Markovian nonlinear dynamical systems subjected to FGN excitations.}

The memFPK equation includes time-dependent exponential integral conditional expectation terms, which present substantial challenges for direct analytical solution. The Volterra adjustable decoupling approximation (VADA) proposed by Mamis et al. \cite{mamis2019systematic} can be used to tackle this difficulty. Here, the exponential integral function is decomposed into its mean and fluctuation (around the mean) components. By then conducting a double expansion on the fluctuation terms and truncating higher-order components, the memory kernel dependence terms are reconstructed in a state-dependent form that removes their reliance on historical states. This yields an approximate memFPK equation, which can be efficiently solved using finite difference (FD) methods to derive the transient response PDF. Although the approximation unavoidably sacrifices some historical details, comparative analyses between numerical experiments and Monte Carlo simulations (MCS) show that the memFPK equation approach keeps a high level of accuracy.  

Overall, a new framework is established for quantifying uncertainty and its propagation in non-Markovian stochastic nonlinear dynamical systems via a novel memFPK equation. The rest of this paper is structured as follows. In Section 2, the main results are introduced. Section 3 offers a detailed derivation of the memFPK equation for nonlinear dynamical systems driven by FGN and GWN is derived via fractional It\^o formula and integration by parts. In Section 4, several numerical examples are presented to verify the effectiveness of the proposed approach.

\section{Problem setting and main result}
Let us consider a nonlinear dynamical system under stochastic excitations of both FGN and GWN. The governing equation of the system is expressed in the following form:
\begin{equation}\label{mSDS}
\dot{X_t}= f(X_t)+g(X_t)\xi_t+h(X_t)\xi_t^H,
\end{equation}
where $ X_t $ represents the nonlinear system response process, and a dot over a variable denotes differentiation with respect to time $ t $; $ f(\cdot), g(\cdot), h(\cdot) $ denote real-valued nonlinear functions; $\xi_t $ is a unit GWN excitation process which is a delta-correlated process having autocorrelation function $ \mathbb{E}[\xi_{t}\xi_{t+\tau}]=\delta(\tau) $; $ \mathbb{E}[\cdot] $ represents the expected value, $ \delta(\cdot) $ is the Dirac's delta function; $\xi_t^{H}$ is a unit FGN excitation process with a power-law autocorrelation function $$ \mathbb{E}[\xi^H_{t}\xi^H_{t+\tau}]=H(2H-1)|\tau|^{2H-2}+2H|\tau|^{2H-1}\delta(\tau),$$ and $\xi^{H}$ and $\xi_t$ are independent; $\xi_t^{H}$ reduces to be $\xi_t$ when $ H = 1/2 $. It is important to note that this paper solely focuses on the one-dimensional case.
		
It is pointed out that when $ 0<H<1/2 $, $\xi^{H}_t$ is not proper for modeling physical noise since its correlation function is negative. In the present paper, only $ 1/2<H<1 $ is considered. In the case of $ 1/2<H<1 $, the power spectral density for $\xi^{H}_t$ is $ S_H(\omega)=\frac{H\Gamma(2H)\sin(H\pi)}{\pi}|\omega|^{1-2H} $, where $ \Gamma(\cdot) $ is Gamma function, for $ H = 1/2 $, the power spectral density reduces to constant power spectrum value equal to $1/(2\pi)$. Thus, we know that $\xi^{H}_t$ has a power-law spectral density function over all frequencies $ \omega $ and exhibits the long-range correlations and non-Markovian characteristics.

The following theorem formulates the memFPK equation corresponding to the  Eq. (\ref{mSDS}). The proof of this theorem will be provided in Section \ref{sec3}.
\begin{theorem}\label{them-1}
Considering Eq. (\ref{mixSDE}) with the initial condition $ X_0=x_0 $ which is a constant or random variable independent of FGN and GWN with known distribution. Suppose that the function  $ f(\cdot) $ is differentiable, while the functions  $ g(\cdot) $ and  $ h(\cdot) $ are twice-differentiable. Additionally, assume that  $ h(x) \neq 0 $  for all  $ x \in R $. Then, the PDF $ p(x,t) $ associated with Eq. (\ref{mixSDE}) is obtained as the solution to the following memFPK equation
\begin{align}\label{FPKnonlinear}
\frac{\partial }{\partial t}p(x,t)=&-\frac{\partial }{\partial x}\big\{a^{{\rm (mem)}}(x,t)p(x,t)\big\}+\frac{\partial^2 }{\partial x^2}\big\{b^{{\rm (mem)}}(x,t)p(x,t)\big\},
\end{align}
where
\begin{align}\label{FPKnonlinear-1}
a^{{\rm (mem)}}(x,t):=&f(x)+\frac{1}{2}g(x)g'(x)+h(x)h'(x)\Psi(x,t),\cr
b^{{\rm (mem)}}(x,t):=&\frac{1}{2}g^2(x)+h^2(x)\Psi(x,t),
\end{align}
with 
\begin{align}\label{FPKnonlinear-2}
\Psi(x,t):=& \mathbb{E}\Big[\int_{0}^{t}\exp\big\{\int_{s}^{t}\varphi_1(X_u)\mathrm{d}u+\int_{s}^{t}\varphi_2(X_u)\mathrm{d}W_u\big\}  \phi(t,s)\mathrm{d}s \mid X_t=x\Big],\cr
\varphi_1:= & f'-f\frac{h'}{h}+\frac{1}{2}g\big(g''-\frac{g'hh'+ghh''-g'(h')^2}{h^2}\big), \quad \varphi_2:= g'-g\frac{h'}{h},
\end{align}
where $\mathrm{d}W_t$ denotes the It\^o integral with respect to Brownian motion (BM) $ W_t $, and $\phi(t,s)=H(2H-1)|t-s|^{2H-2}$. We refer to $a^{{\rm (mem)}}(\cdot,\cdot)$ as the memory-dependent drift coefficient and $b^{{\rm (mem)}}(\cdot,\cdot)$ as the memory-dependent diffusion coefficient. 

Subsequently, we will consider four particular cases that emerge when applying the aforementioned memFPK equation.
\begin{enumerate}[(I)]
\item Provided that $g(x)h'(x) = h(x)g'(x)$ holds true, the memory-dependent drift and diffusion coefficients can be formulated as
\begin{align*}
a^{{\rm (mem)}}(x,t):=&f(x)+\frac{1}{2}g(x)g'(x)\cr
&+h(x)h'(x)\mathbb{E}\Big[\int_{0}^{t}\exp\big\{\int_{s}^{t}\big(f'(X_u)-\frac{h'(X_u)}{h(X_u)}f(X_u)\big)\mathrm{d}u\big\}\phi(t,s)\mathrm{d}s\mid X_t=x\Big],\cr
b^{{\rm (mem)}}(x,t):=&\frac{1}{2}g^2(x)
+h^2(x)\mathbb{E}\Big[\int_{0}^{t}\exp\big\{\int_{s}^{t}\big(f'(X_u)-\frac{h'(X_u)}{h(X_u)}f(X_u)\big)\mathrm{d}u\big\}\phi(t,s)\mathrm{d}s\mid X_t=x\big].
\end{align*}	
\item  Provided that $ f(x)h'(x) = h(x)f'(x) $ and $g(x)h'(x) = h(x)g'(x)$ hold true, the memory-dependent drift and diffusion coefficients can be formulated as
\begin{align*}
a^{{\rm (mem)}}(x,t):=&f(x)+\frac{1}{2}g(x)g'(x)+Ht^{2H-1}h(x)h'(x),\cr
b^{{\rm (mem)}}(x,t):=&\frac{1}{2}g^2(x)+Ht^{2H-1}h^2(x).
\end{align*}	
\item  Provided that $g(x)\equiv\sigma_W, h(x)\equiv\sigma_B$ hold true, the memory-dependent drift and diffusion coefficients can be formulated as
\begin{align*}
a^{{\rm (mem)}}(x,t):=f(x), \quad b^{{\rm (mem)}}(x,t):= \frac{1}{2}\sigma_W^2+\sigma_B^2\mathbb{E}\Big[\int_{0}^{t}\exp\big\{\int_{s}^{t}f'(X_u)\mathrm{d}u\big\}\phi(t,s)\mathrm{d}s\mid X_t=x\Big].
\end{align*}	
\item  Provided that $ f(x) = -\alpha x, g(x)\equiv\sigma_W, h(x)\equiv\sigma_B$ hold true, the memory-dependent drift and diffusion coefficients can be formulated as
\begin{align*}
a^{{\rm (mem)}}(x,t):=-\alpha x, \quad b^{{\rm (mem)}}(x,t):= \frac{1}{2}\sigma_W^2+\sigma_B^2\int_{0}^{t}e^{-\alpha(t-s)}\phi(t,s)\mathrm{d}s,
\end{align*}	
where
\begin{align*}
\int_{0}^{t}e^{-\alpha(t-s)}\phi(t,s)\mathrm{d}s=H(2H-1)\alpha^{1-2H} (\Gamma(2H-1) - \Gamma(2H-1,\alpha t)),
\end{align*}
with $ \Gamma(\alpha,z)=\int_z^{\infty}t^{\alpha-1}e^{-t}\mathrm{d}t $ is upper incomplete Gamma function.
\end{enumerate}
\end{theorem}

\section{Derivation for the memFPK equation}\label{sec3}
\subsection{Overview of FPK equation based on It\^o stochastic calculus theory}
Generally speaking, when  $H=1/2$, the response process 
$X_t$ governed by (\ref{mSDS}) exhibits Markovian behavior. Under this condition, the FPK equation, which serves as the governing equation for the evolution of PDF of $X_t$, can be deduced from Chapman-Kolmogorov equation and Kramers-Moyal expansion, as referenced in \cite{risken1996fokker}.   Alternatively, there is another method to derive the FPK equation using the theoretical framework of It\^o stochastic calculus. Then, the corresponding FPK equation can be derived through It\^o's lemma and integration by parts. The reason for applying It\^o stochastic calculus in this context lies in the fact that, only within the It\^o framework,  the expectation of the stochastic term is zero.  This feature differs sharply from Stratonovich stochastic calculus, where this property does not exist. Importantly, this distinct feature of It\^o stochastic calculus is crucial for deriving the corresponding FPK equation.

Now, we recall a concise overview of the background regarding the derivation of FPK equation for stochastic differential equation (SDE) excited by BM within the theoretical framework of It\^o stochastic calculus. Let us consider the general category of Stratonovich SDE presented in the following form:
\begin{equation}\label{mixSDE-sito}
\mathrm{d}X_t=f(X_t)\mathrm{d}t+g(X_t)\mathrm{d}^{\circ}W_t.
\end{equation}
	
To study the above Stratonovich SDE in the It\^o stochastic calculus framework and derive the corresponding FPK equation, we first introduce the Wong-Zakai  correction. Suppose that  $ X_t $ satisfies the Stratonovich SDE given by (\ref{mixSDE-sito}). For any differentiable function $g(\cdot)$, the distinction between the It\^o integral and the Stratonovich integral is as follows:
\begin{equation}\label{Wong-Zakai-correction}
\int_{0}^{t} g(X_s)  \mathrm{d}^{\circ} W_s= \int_{0}^{t} g(X_s) \mathrm{d}W_s+\frac{1}{2}\int_{0}^{t}g(X_s)g'(X_s)\mathrm{d}s,
\end{equation}
where $\mathrm{d}W_t$ denotes the It\^o integral with respect to $ W_t $, which ensures that the expectation $ \mathbb{E}\big[\int_{0}^{t} g(X_s) \mathrm{d}W_s\big]=0 $ holds. Then, utilizing Eq. (\ref{Wong-Zakai-correction}), we can transform Eq. (\ref{mixSDE-sito}) into the following It\^o SDE:
\begin{equation}\label{mixSDE-ito}
\mathrm{d}X_t=\big(f(X_t)+\frac{1}{2}g(X_t)g'(X_t)\big)\mathrm{d}t+g(X_t)\mathrm{d}W_t,
\end{equation}
where $ \frac{1}{2}\int_{0}^{t}g(X_s)g'(X_s)\mathrm{d}s $ is called the Wong-Zakai correction.

Next, we introduce the It\^o formula pertinent to It\^o SDE. This formula plays a pivotal role in the subsequent derivation of the corresponding FPK equation. Let's assume that the stochastic process $ X_t $ follows the It\^o SDE given in (\ref{mixSDE-ito}). For any function $ F(\cdot) $ that can be differentiated twice, we have
\begin{align}\label{itolemma-bm-sto}
\mathrm{d}F(X_t)=&\big\{\big(f(X_t)+\frac{1}{2}g(X_t)g'(X_t)\big)F'(X_t)+\frac{1}{2}g^2(X_t)F''(X_t)\big\}\mathrm{d}t+g(X_t)F'(X_t)\mathrm{d}W_t.
\end{align}

Upon taking the expectation of Eq. (\ref{itolemma-bm-sto}), we obtain
\begin{align}\label{exp-itolemma-bm}
\frac{\mathrm{d}}{\mathrm{d}t}\mathbb{E}[F(X_t)]=&\mathbb{E}\big[\big(f(X_t)+\frac{1}{2}g(X_t)g'(X_t)\big)F'(X_t)+\frac{1}{2}g^2(X_t)F''(X_t)\big],
\end{align}
which is equivalent to
\begin{align}\label{exp-itolemma-bm-1}
\frac{\mathrm{d} }{\mathrm{d} t}\int_{\mathbb{R}}F(x)p(x,t)\mathrm{d}x=\int_{\mathbb{R}}\bigl\{\big(f(x)+\frac{1}{2}g(x)g'(x)\big)F'(x)+\frac{1}{2}g^2(x)F''(x)\bigr\}p(x,t)\mathrm{d}x.
\end{align}

By performing integration by parts on the right-hand side of Eq. (\ref{exp-itolemma-bm-1}) and assuming that as $ |x|\rightarrow \infty $, $ p(x,t) $ and its derivatives up to the second order with respect to $ x $ all tend to zero, Eq. (\ref{exp-itolemma-bm-1}) can be rewritten as
\begin{align}\label{exp-itolemma-bm-2}	
\int_{\mathbb{R}}F(x)\frac{\partial }{\partial t}p(x,t)\mathrm{d}x
=&
-\int_{\mathbb{R}}F(x) \frac{\partial }{\partial x}\big\{\big(f(x)+\frac{1}{2}g(x)g'(x)\big)p(x,t)\big\}\mathrm{d}x\cr
&+\int_{\mathbb{R}}F(x)\frac{\partial^2 }{\partial x^2} \big\{\frac{1}{2}g^2(x)p(x,t)\big\}\mathrm{d}x.
\end{align}	
	
Making use of the fact that $ F(x) $ is arbitrary, Eq. (\ref{exp-itolemma-bm-2}) leads to the FPK equation
\begin{align}\label{fpk1}	
\frac{\partial }{\partial t}p(x,t)
=
-\frac{\partial }{\partial x}\big\{\big(f(x)+\frac{1}{2}g(x)g'(x)\big)p(x,t)\big\}+\frac{\partial^2 }{\partial x^2} \big\{\frac{1}{2}g^2(x)p(x,t)\big\}.
\end{align}

\subsection{The memFPK equation based on FWIS integral theory}
It is crucial to emphasize that the two derivation methods outlined in Section 3.1 are only applicable to generic Markov processes. However, when $1/2<H<1$, the process  $X_t$ transforms into a non-Markovian process. In this situation, both the Chapman-Kolmogorov equation and the Kramers-Moyal expansion break down. As a consequence, deriving the governing equation for the PDF of FGN-excited system becomes significantly more challenging. At this point, a natural thought arises: could there exists an It\^o-type stochastic calculus theoretical framework tailored to FGN, analogous to the well-established It\^o's stochastic calculus theory for GWN? Such a framework would ideally incorporate corresponding correction terms, an It\^o's formula specific to FGN, and retain the key property that the expected value of the stochastic term is zero. If so, this theoretical construct might potentially enable the derivation of the governing equation that describes the evolution of the PDF for nonlinear dynamical system excited by FGN.

Notably, we make an important discovery, that is, by using FWIS theory \cite{biagini2008stochastic,mishura2008stochastic,ducan2000stochastic}, we have found it possible to derive the Malliavin correction and the fractional It\^o formula. These are similar to the well-known classical ones and are specially designed for the FGN situation.
This breakthrough is not only a turning point but also provides a good basis for great progress in this field. It turns a previously difficult problem into an easy-to-solve one, creating new opportunities for research. In this section, we'll be the first to show the results about the memFPK equation. These results are expected to help us better understand non-Markovian stochastic systems. To make this new contribution more believable, a thorough and careful proof will be given in this section. This proof will be a strong theoretical support, making sure that the proposed ideas and findings are valid and reliable.

The process of deriving the FPK equation for nonlinear dynamical system excited by  a combination of GWN and FGN unfolds in two distinct steps. First, it is necessary to develop a It\^o-type stochastic calculus theory tailored to FGN. This involves exploring properties such as the existence of Wang-Zakai correction, formulating an appropriate version of It\^o lemma, and identifying a stochastic term whose expectation is expected to be 0. Second, employing these newly established properties in conjunction with the integration by parts technique, we can successfully derive the FPK equation. The FWIS integral developed by Duncan and Hu \cite{ducan2000stochastic}, possesses several key properties. It satisfies the condition that the expectation value of the FGN-integral term is 0. Additionally, it adheres to the fractional It\^o lemma and accounts for the difference (Mallivian correction) between the symmetric pathwise integral (Stratonovich type) and the FWIS integral. These properties collectively enable us to derive the equation that describes the evolution of the PDF for nonlinear dynamical system excited by both FGN and GWN. The subsequent discussion will provide a detailed step-by-step derivation of this equation.

We first transform the nonlinear dynamical system excited by FGN and GWN defined by Eq. (\ref{mSDS}) into the following Stratonovich-type SDE:
\begin{equation}\label{mixSDE}
\mathrm{d}X_t=f(X_t) \mathrm{d}t+g(X_t) \mathrm{d}^{\circ}W_t+h(X_t) \mathrm{d}^{\circ} B^H_t,
\end{equation}
where $ W_t$ is a standard BM with $\mathbb{E}[W_t]=0$ and $\mathbb{E}[W_t^2]=t$; $\mathrm{d}^{\circ}W_t$ denotes the Stratonovich integral with respect to the $ W_t $; $ B^H_t$ is a standard fractional Brownian motion (FBM) with $\mathbb{E}[B^H_t]=0$ and $\mathbb{E}[W_t^2]=t^{2H}$; $\mathrm{d}^{\circ}B^H_t$ denotes symmetric pathwise integral \cite{nualart2003stochastic,mishura2008stochastic,biagini2008stochastic} with respect to the $ B_t^H $. $ W_t$ and $ B^H_t$ are independent.

Now, we introduce the correction term in the case of FBM.  Support $ X_t $ satisfies the Stratonovich-type SDE (\ref{mixSDE}), for any differentiable function $ h(\cdot) $, such that 
$$ \mathbb{E}\Big[\int_0^T\int_0^T h(X_t) h(X_s) \phi(t,s)\mathrm{d}s\mathrm{d}t\Big]<\infty,$$ 
with $\phi(t,s)=H(2H-1)|t-s|^{2H-2}$, the differences between the symmetric pathwise integral and FWIS integral is:
\begin{equation}\label{symtofwick}
\int_{0}^{t}h(X_s)\mathrm{d}^{\circ}B^H_s=\int_{0}^{t}h(X_s)\mathrm{d}^{\diamond}B^H_s+\int_{0}^{t}h'(X_s)D^{\phi}_sX_s\mathrm{d}s,
\end{equation}
where $ D^{\phi}_t X_t $ is the Mallivian derivative \cite{ducan2000stochastic} of $X_t$; $\mathrm{d}^{\diamond}B^H$ denotes the FWIS integral with respect to the $ B_t^H$, which ensures that the expectation $ \mathbb{E}\big[\int_{0}^{t}h(X_s)\mathrm{d}^{\diamond}B^H_s\big]=0 $ holds; the term $ \int_{0}^{t}h'(X_s)D^{\phi}_sX_s\mathrm{d}s $ is called the Mallivian correction.

By employing Wong-Zakai and Mallivian corrections, we are able to reformulate the SDE (\ref{mixSDE}) as
\begin{equation}\label{mixSDEfito}
\mathrm{d}X_t=\big(f(X_t)+\frac{1}{2}g(X_t)g'(X_t)+h'(X_t)D^{\phi}_tX_t \big)\mathrm{d}t+g(X_t)\mathrm{d}W_t+h(X_t) \mathrm{d}^{\diamond}B^{H}_t.
\end{equation}	

Next, we introduce the fractional It\^o's lemma (see, Lemma \ref{itomixed}) pertinent to (\ref{mixSDEfito}). Let the stochastic process $ X_t $ satisfy the It\^o-type SDE given by \eqref{mixSDEfito}. For any function $ F(\cdot) $ that is twice differentiable, one has 
\begin{align}\label{itolemma-fbm}
\mathrm{d}F(X_t)=&\big\{\big(f(X_t)+\frac{1}{2}g(X_t)g'(X_t)+h'(X_t)D^{\phi}_tX_t\big)F'(X_t)\cr
&+\big(\frac{1}{2}g^2(X_t)+h(X_t)D^{\phi}_tX_t\big)F''(X_t)\big\}\mathrm{d}t+g(X_t)F'(X_t)\mathrm{d}W_t+h(X_t)F'(X_t) \mathrm{d}^{\diamond}B^H_t.
\end{align}
{\bf The Proof of Theorem \ref{them-1}:} 
We begin by considering the case of constant FGN diffusion, where $h(x)\equiv\sigma_B$. In this case, the SDE given by \eqref{mixSDE} can be simplified to the following form:
\begin{equation}\label{mixSDE-add}
\mathrm{d}X_t= f(X_t) \mathrm{d}t+ g(X_t) \mathrm{d}^{\circ}W_t+\sigma_B \mathrm{d}^{\circ}B^{H}_t.
\end{equation}	

Consequently, it is straightforward to derive its It\^o-type SDE along with the fractional It\^o formula, i.e.
\begin{equation}\label{mixSDE-add-ito}
\mathrm{d}X_t=\big(f(X_t)+\frac{1}{2}g(X_t)g'(X_t)\big)\mathrm{d}t+g(X_t)\mathrm{d}W_t+\sigma_B \mathrm{d}^{\diamond}B^{H}_t,
\end{equation}
and
\begin{align}\label{itolemma-fbmadd}
\mathrm{d}F(X_t)=&\big\{\big(f(X_t)+\frac{1}{2}g(X_t)g'(X_t)\big)F'(X_t)+\big(\frac{1}{2}g^2(X_t)+\sigma_BD^{\phi}_tX_t\big)F''(X_t)\big\}\mathrm{d}t\cr
&+g(X_t)F'(X_t)\mathrm{d}W_t+\sigma_B F'(X_t) \mathrm{d}^{\diamond}B^H_t.
\end{align}

Taking the expectation on Eq. (\ref{itolemma-fbmadd}), one has 
\begin{align}\label{fp1-1}
\frac{\mathrm{d}}{\mathrm{d}t}\mathbb{E}[F(X_t)]=&
\mathbb{E}\big[\big(f(X_t)+\frac{1}{2}g(X_t)g'(X_t)\big)F'(X_t)\big]+\mathbb{E}\big[\big(\frac{1}{2}g^2(X_t)+\sigma_B D^{\phi}_tX_t\big)F''(X_t)\big],
\end{align}
where the explicit formulation of $D^{\phi}_t X_t$ is
\begin{align}\label{Malfbmx-1}
D^{\phi}_tX_t=
\sigma_B\int_{0}^{t}\exp\big\{\int_{s}^{t}\big(f'(X_u)-\frac{1}{2}g(X_u)g''(X_u)\big)\mathrm{d}u+\int_{s}^{t}g'(X_u)\mathrm{d}W_u\big\} \phi(t,s)\mathrm{d}s.
\end{align}

It is evident that $ D^{\phi}_tX_t $ exhibits history-dependence, thus, we called it as a memory kernel function. For a detailed derivation of Eq. (\ref{Malfbmx-1}), the interested reader is encouraged to refer to \ref{app-a1}.

\begin{remark}\label{mal-3cases}
Let $ X_t $ satisfy the Stratonovich-type SDE given by \eqref{mixSDE-add}. Three special cases are noted here.
\begin{itemize}
        	
\item If $g(x)\equiv\sigma_W, h(x)\equiv\sigma_B$ hold, then the Mallivian derivative becomes
\begin{equation}\label{Mal-case1}
D^{\phi}_tX_t= \sigma_B\int_{0}^{t}e^{\int_{s}^{t}f'(X_u)\mathrm{d}u}\phi(t,s)\mathrm{d}s.
\end{equation}
        	
\item If $ f(x) = -\alpha x, g(x)\equiv\sigma_W, h(x)\equiv\sigma_B$ hold, then the Mallivian derivative becomes
\begin{equation}\label{Mal-case2}
D^{\phi}_tX_t= \sigma_B\int_{0}^{t}e^{-\alpha (t-s)}\phi(t,s)\mathrm{d}s.
\end{equation}

\item If $ f(x) = -\alpha, g(x)\equiv\sigma_W, h(x)\equiv\sigma_B$ hold, then the Mallivian derivative becomes
\begin{equation}\label{Mal-case3}
	D^{\phi}_tX_t= \sigma_B Ht^{2H-1}.
\end{equation}

\end{itemize}
\end{remark}

By applying the integration by parts technique once more, Eq. (\ref{fp1-1}) can be derived in a manner analogous to the BM case, as demonstrated in Eqs. (\ref{exp-itolemma-bm-2}) and (\ref{fpk1}). This process results in the memFPK equation as follows:
\begin{align}\label{FPKnonlinear-add}
\frac{\partial }{\partial t}p(x,t)=-\frac{\partial }{\partial x}\big\{a^{{\rm (mem)}}(x,t)p(x,t)\big\}+\frac{\partial^2 }{\partial x^2}\big\{b^{{\rm (mem)}}(x,t)p(x,t)\big\},
\end{align}
where the memory-dependent drift and diffusion coefficient are expressed as
\begin{align}\label{FPKnonlinear-add-1}
	a^{{\rm (mem)}}(x,t)&:=f(x)+\frac{1}{2}g(x)g'(x),\cr
	b^{{\rm (mem)}}(x,t)&:=\frac{1}{2}g^2(x)+\sigma_B^2 \Phi(x,t),
\end{align}
with the term relying on the memory kernel
\begin{align}\label{FPKnonlinear-add-2}
\Phi(x,t):=\mathbb{E}\Big[\int_{0}^{t}\exp\big\{\int_{s}^{t}\Big(f'(X_u)-\frac{1}{2}g(X_u)g''(X_u)\Big)\mathrm{d}u+\int_{s}^{t}g'(X_u)\mathrm{d}W_u\big\} \phi(t,s)\mathrm{d}s \mid X_t=x\Big].
\end{align}

Subsequently, we employ the Lamperti transformation to analyze the nonlinear FGN diffusion case, specifically focusing on the SDE given by (\ref{mixSDE}). Through this approach, we are able to derive the memFPK equation. Suppose that $ h(x) \neq 0 $ for all $ x \in \mathbb{R} $. Let us define the function 
$ G(x):=\int_{0}^{x} \frac{1}{h(z)} \mathrm{d}z $, which is known as the Lamperti transform. Since $ h(x) \neq 0 $ for all $ x \in \mathbb{R} $, the function $ G $ is both monotonic and invertible. We denote its inverse function as $ G^{-1} $. A straightforward calculation reveals that
\begin{equation}\label{eq2-2-2}
G'(x)=\frac{1}{h(x)},\qquad (G^{-1})'(x)=\frac{1}{G'(G^{-1}(x))}=h(G^{-1}(x)).
\end{equation}
      
To proceed, we denote
\begin{equation}\label{eq2-2-3}
Y_t:=G(X_t)=\int_{0}^{X_t} \frac{1}{h(z)} \mathrm{d}z,
\end{equation}
and use fractional It\^o's lemma for Stratonovich-type SDE (\ref{mixSDE}) (see, Lemma \ref{stomixed}) to have
\begin{align}\label{eq2-2-4}
\mathrm{d}Y_t&=G'(X_t)\big(f(X_t) \mathrm{d}t+g(X_t) \mathrm{d}^{\circ}W_t+h(X_t) \mathrm{d}^{\circ} B^H_t\big)=\frac{f(X_t)}{h(X_t)} \mathrm{d}t+\frac{g(X_t)}{h(X_t)} \mathrm{d}^{\circ}W_t+  \mathrm{d}^{\circ} B^H_t\cr
&=\frac{f(G^{-1}(Y_t))}{h(G^{-1}(Y_t))} \mathrm{d}t+\frac{g(G^{-1}(Y_t))}{h(G^{-1}(Y_t))} \mathrm{d}^{\circ}W_t+  \mathrm{d}^{\diamond} B^H_t.
\end{align}

For simplicity, we denote
\begin{align}\label{coff-fg}
\tilde{f}(Y_s):=\frac{f(G^{-1}(Y_s))}{h(G^{-1}(Y_s))},\quad \tilde{g}(Y_s):=\frac{g(G^{-1}(Y_s))}{h(G^{-1}(Y_s))}.
\end{align}

Taking into account Eq. (\ref{coff-fg}) into (\ref{eq2-2-4}) yields
\begin{equation}\label{eq2-2-5}
\mathrm{d}Y_t=\tilde{f}(Y_t) \mathrm{d}t+ \tilde{g}(Y_t) \mathrm{d}^{\circ}W_t+\mathrm{d}^{\diamond} B^H_t.
\end{equation}

Building on the derivations from the constant FGN diffusion case, we can obtain the memFPK equation for the stochastic process  $Y_t$
\begin{align}\label{eq2-2-6}
\frac{\partial }{\partial t}p(y,t)=&-\frac{\partial }{\partial y}\big\{a_Y^{{\rm (mem)}}(y,t)p(y,t)\big\}+\frac{\partial^2 }{\partial y^2}\big\{b_Y^{{\rm (mem)}}(y,t)p(y,t)\big\},
\end{align}
where the memory-dependent drift and diffusion coefficients are, respectively,
\begin{align}\label{eq2-2-6-1}
a_Y^{{\rm (mem)}}(y,t)&=\tilde{f}(y)+\frac{1}{2}\tilde{g}(y)\tilde{g}'(y),\cr	
b_Y^{{\rm (mem)}}(y,t)&=\frac{1}{2}\tilde{g}^2(y)+ \mathbb{E}\big[D^{\phi}_t Y_t  \mid Y_t=y\big],
\end{align}
with
\begin{align}\label{eq2-2-7}
D^{\phi}_t Y_t&=\int_{0}^{t}\exp\big\{\int_{s}^{t}\Big(\tilde{f}'(Y_u)-\frac{1}{2}\tilde{g}(Y_u)\tilde{g}''(Y_u)\Big)\mathrm{d}u+\int_{s}^{t}\tilde{g}'(Y_u)\mathrm{d}W_u\big\} \phi(t,s)\mathrm{d}s\cr
&=\int_{0}^{t}\exp\big\{\int_{s}^{t}\tilde{\varphi}_1(Y_u)\mathrm{d}u+\int_{s}^{t}\tilde{\varphi}_2(Y_u)\mathrm{d}W_u\big\} \phi(t,s)\mathrm{d}s,
\end{align}
where
\begin{align*}
\begin{cases}
\tilde{f}'(y)=\big(\big(f'-f\frac{h'}{h}\big)\circ G^{-1}\big)(y),\cr
\tilde{g}'(y)=\big(\big(g'-g\frac{h'}{h}\big)\circ G^{-1}\big)(y),\cr
\tilde{g}''(y)=\big(h\big(g''-\frac{g'hh'+ghh''-g'(h')^2}{h^2}\big)\circ G^{-1}\big)(y),
\end{cases}
\end{align*}
and
\begin{align}\label{eq2-2-8}
\begin{cases}
\tilde{\varphi}_1(y):=\big(\big(f'-f\frac{h'}{h}-\frac{1}{2}g\big(g''-\frac{g'hh'+ghh''-g'(h')^2}{h^2}\big)\big)\circ G^{-1}\big)(y),\cr
\tilde{\varphi}_2(y):=\big(\big(g'-g\frac{h'}{h}\big)\circ G^{-1}\big)(y).
\end{cases}
\end{align}
Here, $ F_1 \circ F_2 $ deontes the composite of two functions $ F_1, F_2 $ of $ \mathbb{R} $ into itself, i.e., $ (F_1 \circ F_2)(\cdot) = F_1(F_2(\cdot))$.

Let $p(x,t)$ denote the PDF for $X_t$ in (\ref{mixSDE}) and the initial condition $ X_0=x_0 $ be a constant or random variable with known distribution. Using the relation $ p(y,t)=h(x)p(x,t) $ and $ y=G(x) $, Eq. (\ref{eq2-2-6}) transforms into
\begin{align}\label{eq2-2-9}
h(x)\frac{\partial }{\partial t}p(x,t)=&-\frac{\partial }{\partial G(x)}\big\{\big(\tilde{f}(G(x))+\frac{1}{2}\tilde{g}(G(x))\tilde{g}'(G(x))\big)h(x)p(x,t)\big\}\cr
&+\frac{\partial^2 }{\partial G^2(x)}\big\{\big(\frac{1}{2}\tilde{g}^2(G(x))+\mathbb{E}\big[D^{\phi}_t G(X_t)  \mid G(X_t)=G(x)\big]\big)h(x)p(x,t)\big\}\cr
=&-h(x)\frac{\partial }{\partial x}\big\{\big(\tilde{f}(G(x))+\frac{1}{2}\tilde{g}(G(x))\tilde{g}'(G(x))\big)h(x)p(x,t)\big\}\cr
&+h(x)\frac{\partial }{\partial x}\big\{h(x)\frac{\partial }{\partial x}\big\{\big(\frac{1}{2}\tilde{g}^2(G(x))+\mathbb{E}\big[D^{\phi}_t G(X_t)  \mid X_t=x\big]\big)h(x)p(x,t)\big\}\big\},
\end{align}
where the last equality relies on the assumption that the function $ G $ is both monotonic and invertible, and that  $ X_t $ is a real-valued stochastic process. Under these conditions, the events $ \{G(X_t)=G(x)\} $ and $ \{X_t=x\} $ are equivalent, which ensures that the equality $ \mathbb{E}\big[D^{\phi}_t G(X_t) \mid G(X_t)=G(x)\big]=\mathbb{E}\big[D^{\phi}_t G(X_t)  \mid X_t=x\big] $ holds.

Next, integrating by parts the terms in Eq. (\ref{eq2-2-9}) yields 
\begin{align}\label{eq2-2-10}
\frac{\partial }{\partial t}p(x,t)=&-\frac{\partial }{\partial x}\big\{\big(\tilde{f}(G(x))+\frac{1}{2}\tilde{g}(G(x))\tilde{g}'(G(x))\big)h(x)p(x,t)\big\}\cr
&+\frac{\partial }{\partial x}\big\{h(x)\frac{\partial }{\partial x}\big\{\big(\frac{1}{2}\tilde{g}^2(G(x))+\mathbb{E}\big[D^{\phi}_t G(X_t)  \mid X_t=x\big]\big)h(x)p(x,t)\big\}\big\}\cr
=&-\frac{\partial }{\partial x}\big\{\big(\tilde{f}(G(x))+\frac{1}{2}\tilde{g}(G(x))\tilde{g}'(G(x))\big) h(x)p(x,t)\big)\big\}\cr
&-\frac{\partial }{\partial x}\big\{h'(x)\big(\frac{1}{2}\tilde{g}^2(G(x))+\mathbb{E}\big[D^{\phi}_t G(X_t)  \mid X_t=x\big]\big)h(x)p(x,t)\big\}\cr
&+\frac{\partial^2 }{\partial x^2}\big\{h(x)\big(\frac{1}{2}\tilde{g}^2(G(x))+\mathbb{E}\big[D^{\phi}_t G(X_t)  \mid X_t=x\big]\big)h(x)p(x,t)\big\}\cr
=&-\frac{\partial }{\partial x}\big\{\big(f(x)+\frac{1}{2}g(x)g'(x)\big)p(x,t)\big\}-\frac{\partial }{\partial x}\big\{h(x)h'(x)\mathbb{E}\big[D^{\phi}_t G(X_t)  \mid X_t=x\big]p(x,t)\big\}\cr
&+\frac{\partial^2 }{\partial x^2}\big\{\big(\frac{1}{2}g^2(x)+h^2(x)\mathbb{E}\big[D^{\phi}_t G(X_t)  \mid X_t=x\big]\big)p(x,t)\big\},
\end{align}
with
\begin{align*}
D^{\phi}_t G(X_t) =\int_{0}^{t}\exp\big\{\int_{s}^{t}\tilde{\varphi}_1(G(X_u))\mathrm{d}u+\int_{s}^{t}\tilde{\varphi}_2(G(X_u))\mathrm{d}W_u\big\} \phi(t,s)\mathrm{d}s,
\end{align*}
where
\begin{align*}
\tilde{\varphi}_1(G(x))&=\big(\tilde{\varphi}_1\circ G\big)(x)\cr
&=\big(\big(f'-f\frac{h'}{h}-\frac{1}{2}g\big(g''-\frac{g'hh'+ghh''-g'(h')^2}{h^2}\big)\big)\circ G^{-1}\circ G\big) (x)\cr
&=\big(f'-f\frac{h'}{h}-\frac{1}{2}g\big(g''-\frac{g'hh'+ghh''-g'(h')^2}{h^2}\big)(x)=:\varphi_1(x),
\end{align*}
and
\begin{align*}
\tilde{\varphi}_2(G(x))=\big(\tilde{\varphi}_2'\circ G\big)(x)
=\big(\big(g'-g\frac{h'}{h}\big)\circ G^{-1}\circ G\big) (x)
=\big(g'-g\frac{h'}{h}\big)(x)
=:\varphi_2(x).
\end{align*}
Up to now,  we obtain the desired memFPK equation for nonlinear SDE (\ref{mixSDE}) as in Theorem \ref{them-1}.
{			
\section{Step-by-step procedure of the proposed method}
In summary, a memFPK equation method for probability response analysis of nonlinear dynamical systems subjected to FGN and GWN is proposed in the present paper, whose schematic flowchart is illustrated in Fig. \ref{flowchart}. The detailed steps outlined below:
\begin{figure}[h]
	\centering
	\subfloat{\includegraphics[scale=0.10]{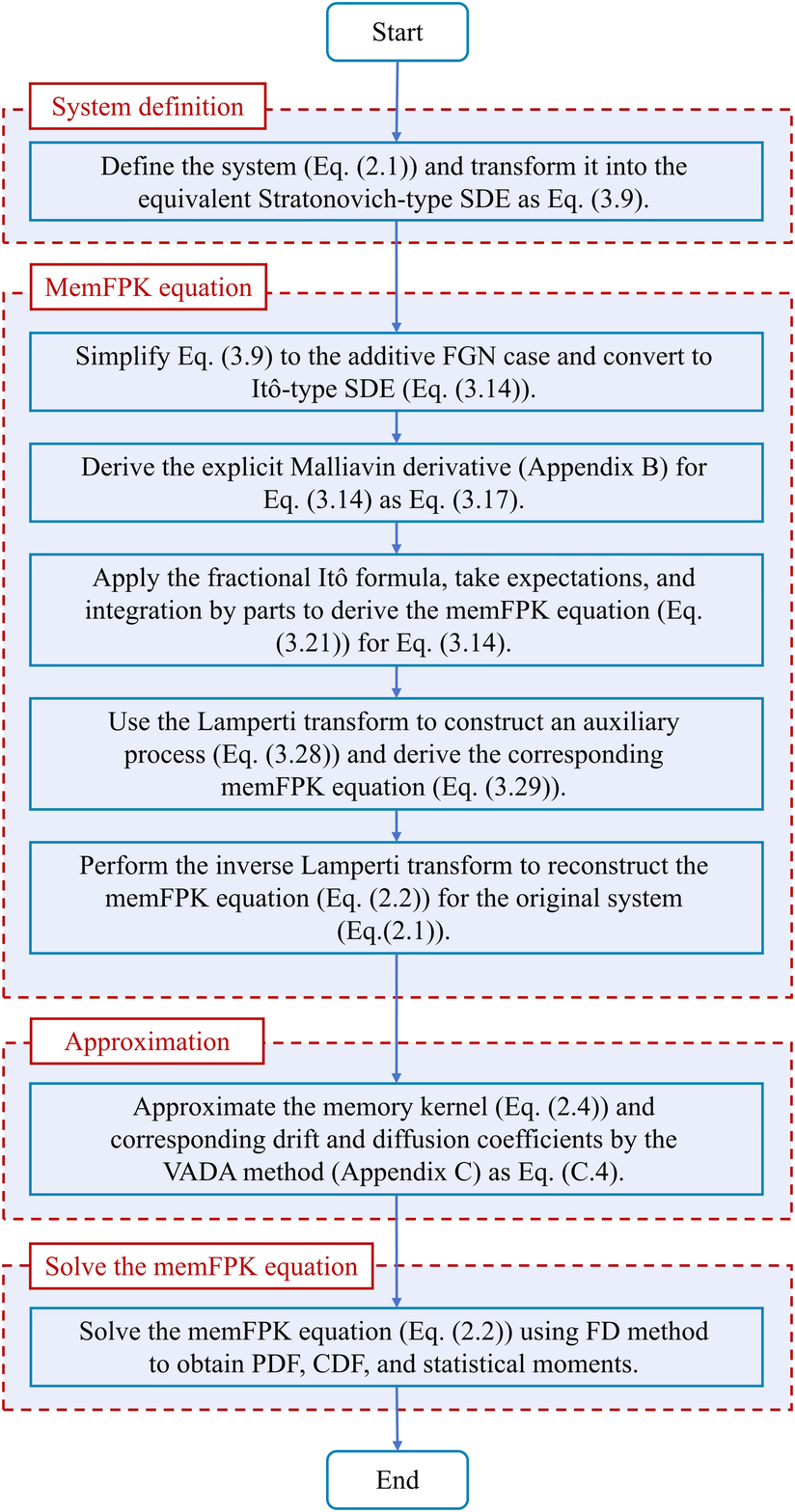}}
	\caption{Flowchart of the proposed method.}
	\label{flowchart}
\end{figure}

\begin{enumerate}[\textit{\textbf{Step}} 1]
	\item Specify the original nonlinear dynamical system as Eq. (\ref{mSDS}) and reformulate it equivalently to Stratonovich-type SDE as Eq.(\ref{mixSDE}).
	\item For tractability, simplify Eq.(\ref{mixSDE}) to the additive FGN case ($h(x)\equiv\sigma_B$) as Eq.(\ref{mixSDE-add}), and convert the Stratonovich-type form to It\^o-type SDE (Eq.(\ref{mixSDE-add-ito})).
	\item Apply the fractional It\^o formula to function $F(X_t)$ to obtain Eq.(\ref{itolemma-fbmadd}), then take expectations as Eq.(\ref{fp1-1}) and derive the explicit formulation (Appendix B) of Malliavin derivative $D_t^{\phi}X_t$ as Eq. (\ref{Malfbmx-1}).
	\item Perform integration by parts on Eq.(\ref{fp1-1}) to obtain the memFPK equation (Eq.(\ref{FPKnonlinear-add})), and the memory-dependent drift and diffusion coefficients (Eq.(\ref{FPKnonlinear-add-1})) together with the memory kernel (Eq.(\ref{FPKnonlinear-add-2})).
	\item Apply a Lamperti transform to the general multiplicative Stratonovich-type SDE (Eq.(\ref{mixSDE})) to construct an auxiliary additive process as Eq.(\ref{eq2-2-5}).
	\item Repeat Steps 2-4 for the auxiliary additive process (Eq.(\ref{eq2-2-5})) to obtain the corresponding memFPK equation (Eq.(\ref{eq2-2-6})) and memory-dependent coefficients (Eq.(\ref{eq2-2-6-1})).
	\item Perform the inverse Lamperti transform and integration by parts to reconstruct the memFPK equation for the original system (Eq. (\ref{mSDS})) as Eq.(\ref{FPKnonlinear}) and its associated memory-dependent coefficients and kernel (Eq.(\ref{FPKnonlinear-1})- (\ref{FPKnonlinear-2})).
	\item In the commutative case of diffusion coefficients (i.e., $g'(x)h(x)=g(x)h'(x)$), apply the VADA method (Appendix C) to the memory kernel (Eq.(\ref{FPKnonlinear-2})), yielding the approximate kernel as Eq.(\ref{eq3-3}) and corresponding drift and diffusion coefficients.
	\item Substitute the approximate kernel (Eq.(\ref{eq3-3})) into the memFPK equation (Eq.(\ref{FPKnonlinear})) and numerically solve it using the FD method to obtain the estimated PDF, CDF and statistical moments.
\end{enumerate} 
}

\section{Numerical results}
In the following section, the effectiveness of the proposed memFPK equation will be verified through numerical examples.

\subsection{Ornstein-Uhlenbeck process}
Let us consider an Ornstein-Uhlenbeck process. Its dynamics are described by the following form: 
\begin{equation}\label{example-1}
\dot{X_t}= -\alpha X_t+\sigma_W\xi_t+\sigma_B\xi_t^H,
\end{equation}
where $ \alpha > 0 $ is a linear coefficient, $ \xi_t $ and $ \xi^H_t $ are unit GWN and FGN respectively, while $ \sigma_W^2 $ and $ \sigma_B^2 $ are their noise intensities respectively. The initial value is taken as $X_0=x_0$, a constant. 
		
Then, refer to Theorem \ref{them-1} (IV), the response PDF of Eq. (\ref{example-1}) satisfies the memFPK equation (\ref{FPKnonlinear}) with the following memory-dependent drift and diffusion coefficients:
\begin{align*}
a^{{\rm (mem)}}(x,t)=-\alpha x, \quad b^{{\rm (mem)}}(x,t)= \frac{1}{2}\sigma_W^2+\sigma_B^2\int_{0}^{t}e^{-\alpha(t-s)}\phi(t,s)\mathrm{d}s.
\end{align*}
	
The analytical solution of the response PDF is given as
\begin{align}\label{fpk-case2-exact-1}
p(x,t)=\frac{1}{\sqrt{2\pi \sigma_x^2(t)}}\exp\Bigl\{-\frac{(x-\mu_x(t))^2}{2\sigma_x^2(t)}\Bigr\},
\end{align}
where
\begin{equation*}
\begin{cases}
\mu_x(t)=&x_0e^{-\alpha t},\\
\sigma_x^2(t)=&\frac{\sigma_W^2}{2\alpha}(1-e^{-2\alpha t})+\sigma_B^2\int_{0}^{t}\int_{0}^{t}e^{-\alpha(t-u)}e^{-\alpha(t-v)}\phi(u,v)\mathrm{d}u\mathrm{d}v.
\end{cases}
\end{equation*}

The memFPK equation can be numerically calculated by FD method. {In the following, the parameters of the Ornstein-Uhlenbeck process are taken as table \ref{table-1}.} The computing domain in solving the memFPK equation is taken as $ [-5,5] $ with the space step $ \Delta x=0.05 $, and the time step $ \Delta t $ takes the value 0.0001 s.
\begin{table}[h!]
	\centering
	{\caption{Values of parameters of Ornstein-Uhlenbeck process for example 5.1}
	\begin{tabular}{ccccc}
		\toprule
		parameters & $\alpha$ & $\sigma_W$ & $\sigma_B$ & $x_0$ \\ 
		\midrule
		values & 1 & 1 & 1 & 1  \\ 
		\bottomrule
	\end{tabular}
	\label{table-1}}
\end{table}

\begin{figure}[h]
\centering
\subfloat[PDF ($ H=0.65 $)]{
	\includegraphics[scale=0.50]{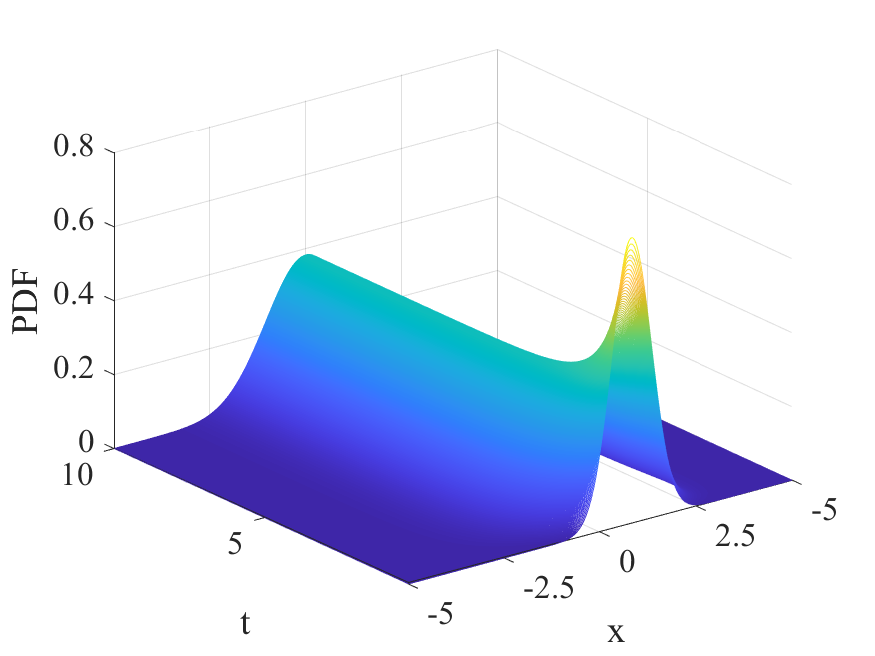}}
\subfloat[Analytical PDF ($ H=0.65 $)]{
	\includegraphics[scale=0.50]{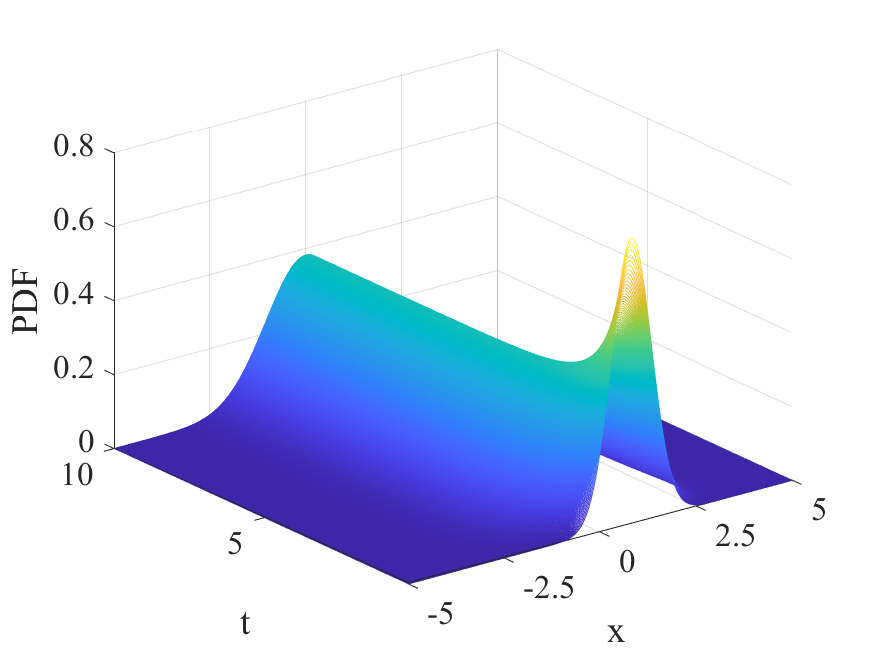}}
\\
\subfloat[Error surface ($ H=0.65 $)]{
	\includegraphics[scale=0.50]{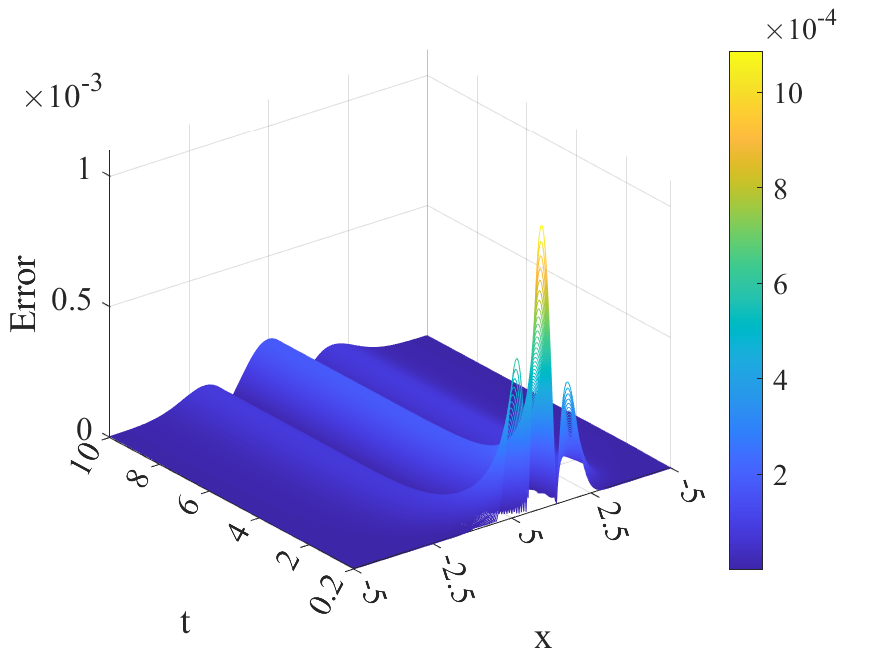}}
\subfloat[Error value($ H=0.65 $)]{
	\includegraphics[scale=0.50]{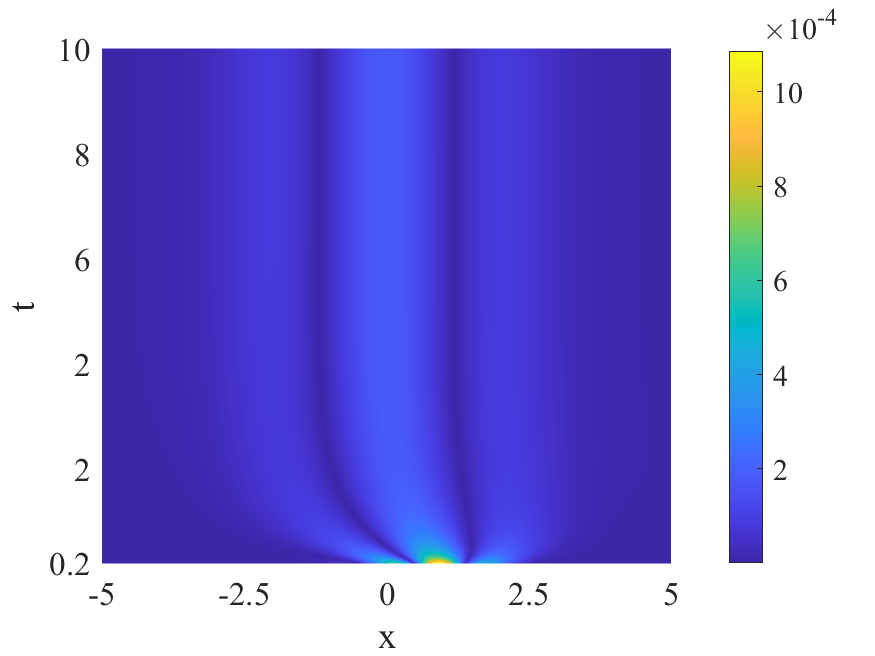}}
\caption{Comparison of the PDF surface via memFPK equation.}
\label{fig3-1-3-1}
\end{figure}

\begin{figure}[t!]
\centering
\subfloat[PDF ($ H=0.95 $)]{
	\includegraphics[scale=0.50]{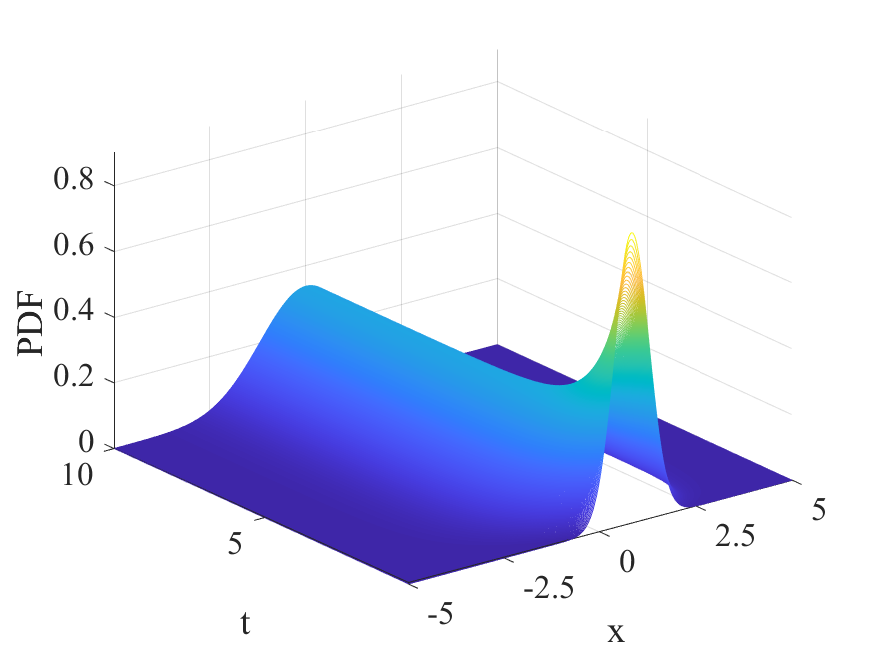}}
\subfloat[Analytical PDF ($ H=0.95 $)]{
	\includegraphics[scale=0.50]{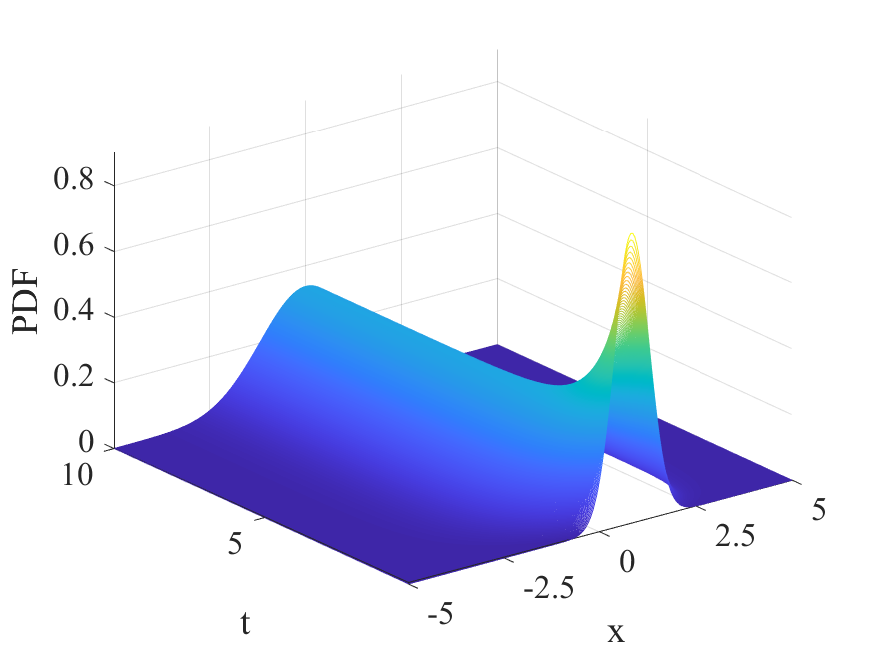}}
\\
\subfloat[Error surface($ H=0.95 $)]{
	\includegraphics[scale=0.50]{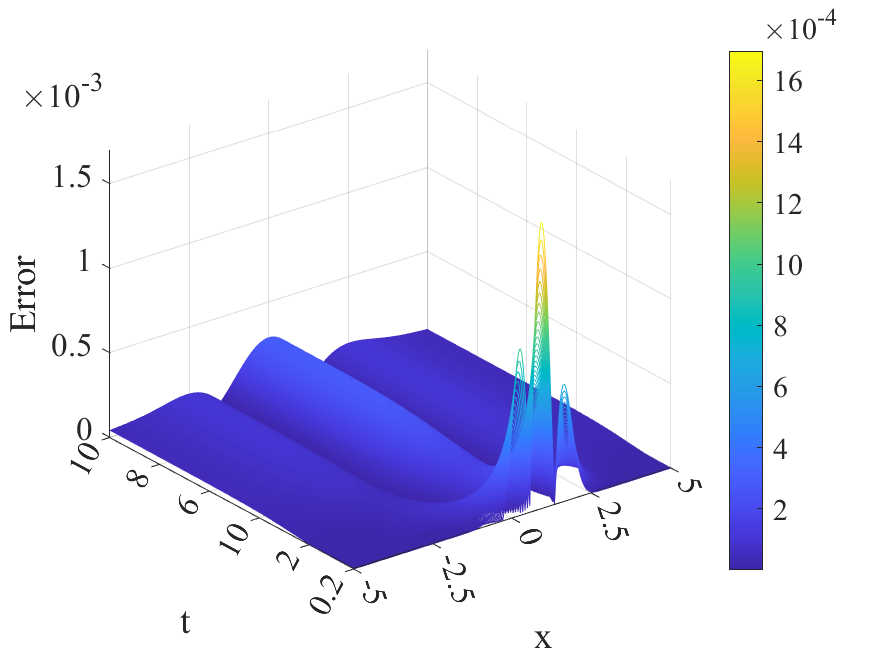}}
\subfloat[Error value ($ H=0.95 $)]{
	\includegraphics[scale=0.50]{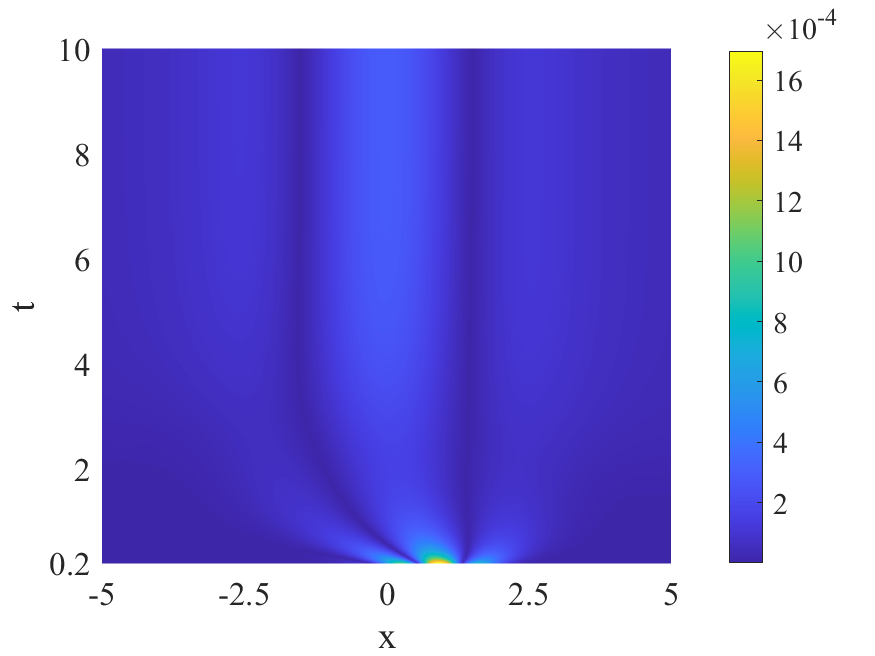}}
\caption{Comparison of the PDF surface via memFPK equation.}
\label{fig3-1-3-2}
\end{figure}

In this example, we examine two distinct Hurst parameters, $H=0.65$ and $H=0.95$. Figs. \ref{fig3-1-3-1} and \ref{fig3-1-3-2} illustrate the evolution surfaces of the PDF for these different $H$ values. In Fig. \ref{fig3-1-3-1}, panel (a) shows the evolution surface of the PDF solution obtained from the memFPK equation, while panel (b) depicts the evolution surface of the analytical PDF solution. Panels (c) and (d) of Fig. \ref{fig3-1-3-1} present the absolute error evolution diagrams between these two methods for 
$H=0.65$. Similarly, for $H=0.95$, the corresponding PDF evolution surfaces and their associated absolute error analysis results are shown in panels (a) and (b), and panels (c) and (d) of Fig. \ref{fig3-1-3-2}, respectively.

\begin{figure}[t!]
\centering
\subfloat[linear scale ($ H=0.65 $)]{
	\includegraphics[scale=0.50]{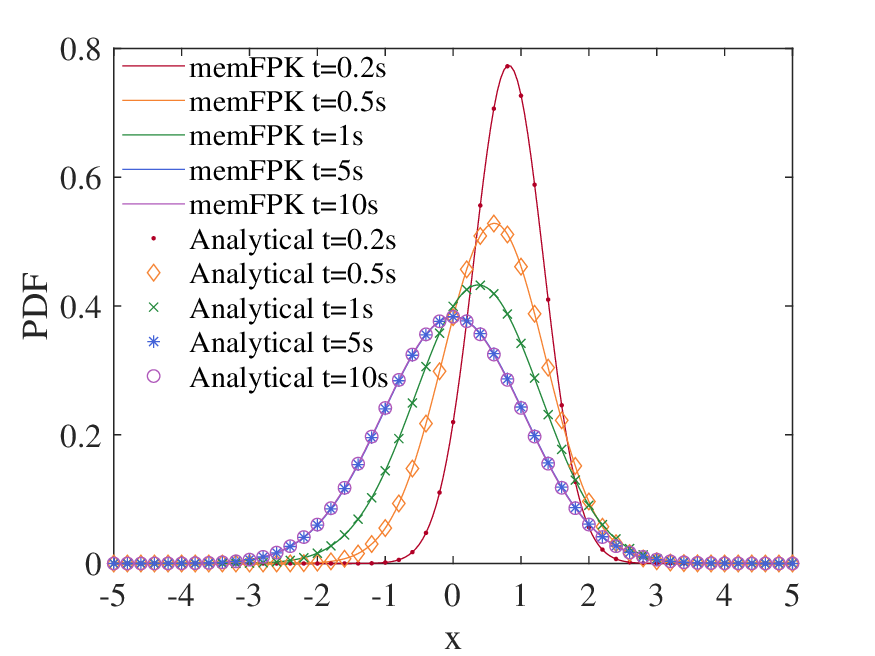}}
\subfloat[logarithmic scale ($ H=0.65 $)]{
	\includegraphics[scale=0.50]{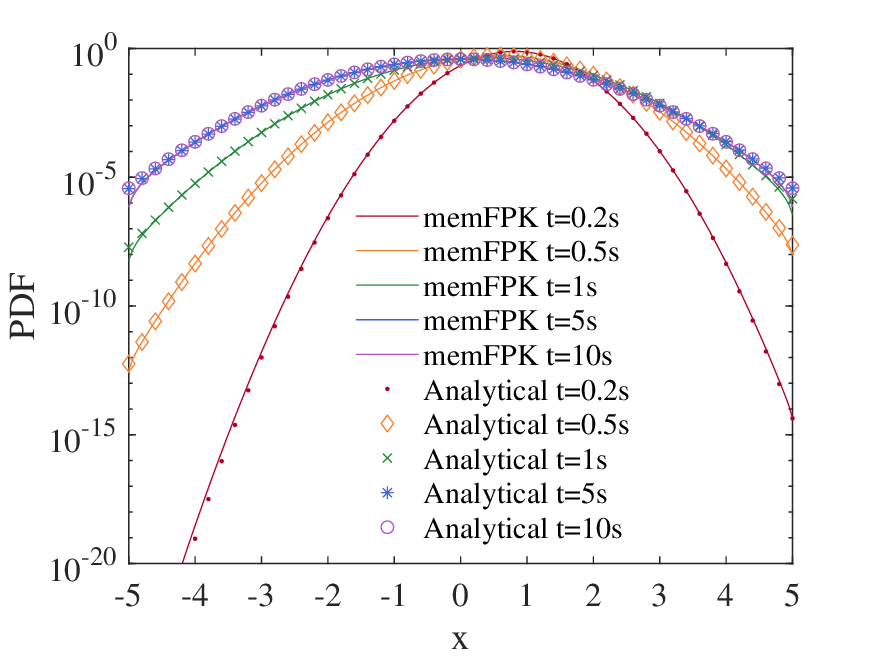}}
\\
\subfloat[linear scale ($ H=0.95 $)]{
	\includegraphics[scale=0.50]{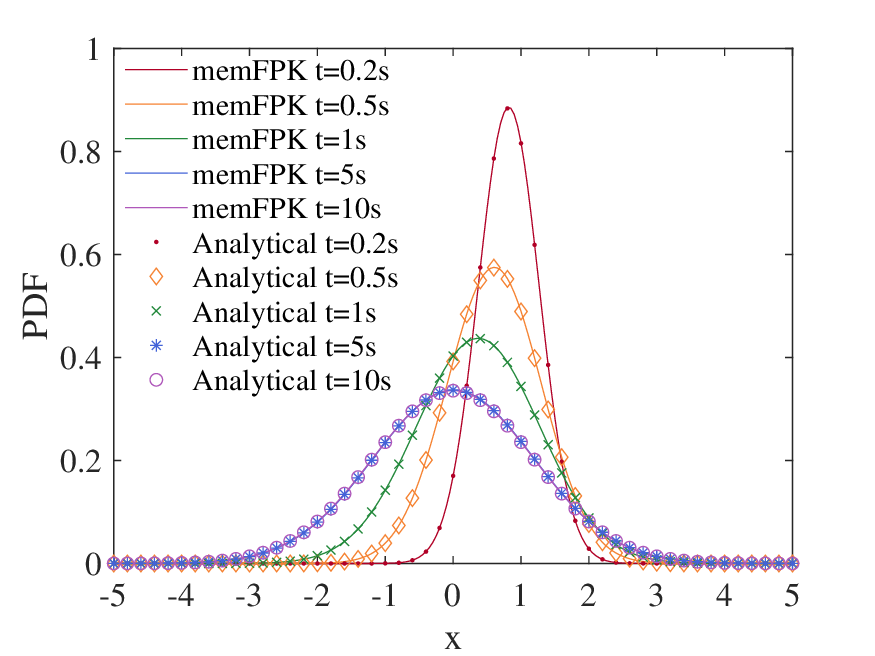}}
\subfloat[logarithmic scale ($ H=0.95 $)]{
	\includegraphics[scale=0.50]{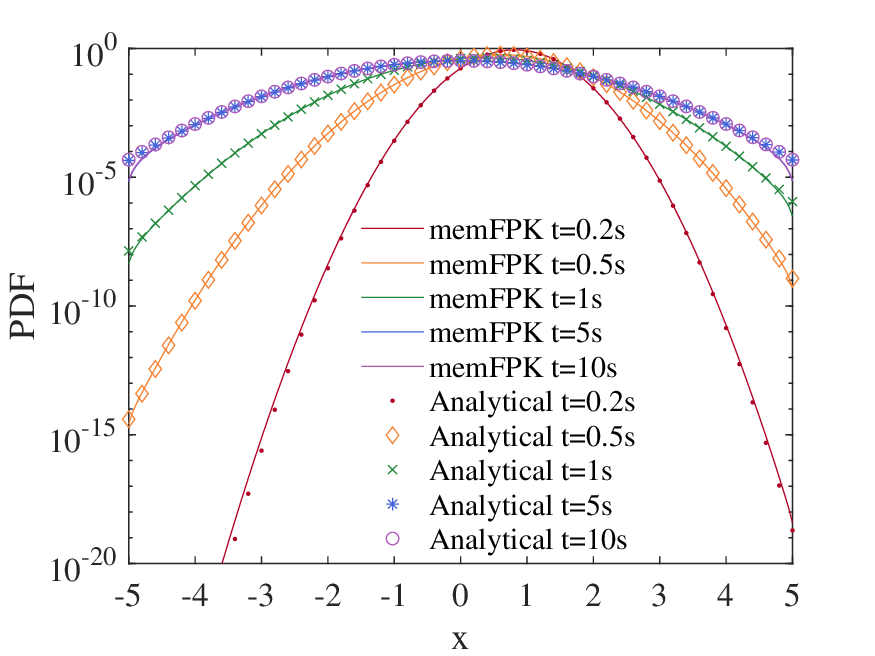}}
\caption{Comparison of the transient PDFs via memFPK equation.}
\label{fig3-1-1}
\end{figure}

\begin{figure}[t!]
\centering
\subfloat[linear scale ($ H=0.65 $)]{
	\includegraphics[scale=0.50]{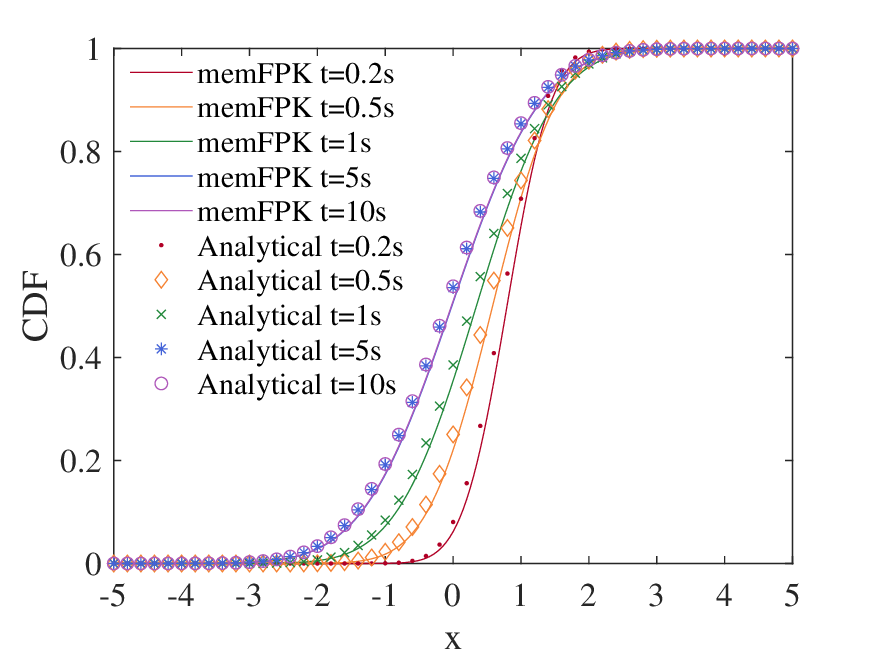}}
\subfloat[logarithmic scale ($ H=0.65 $)]{
	\includegraphics[scale=0.50]{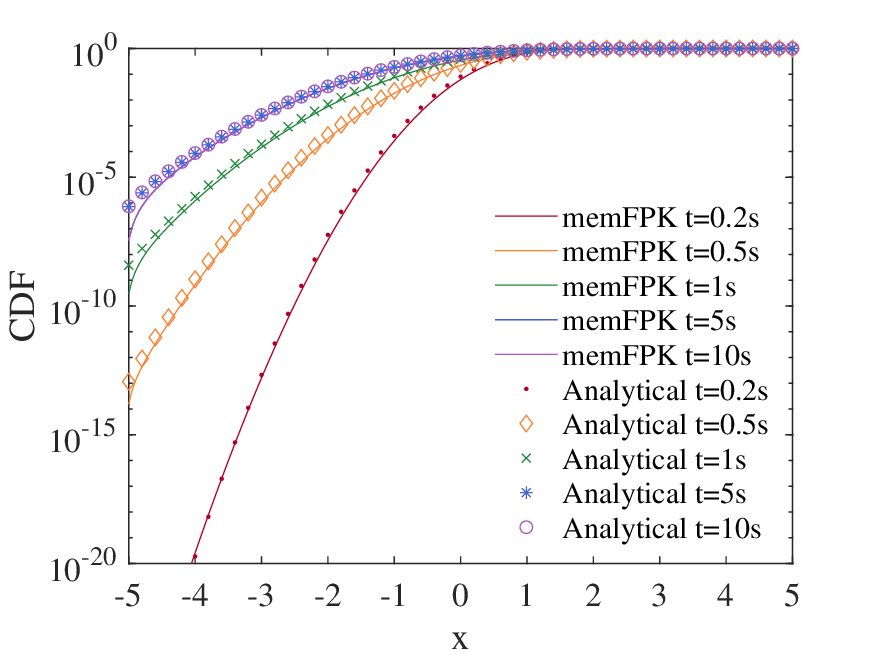}}
\\
\subfloat[linear scale ($ H=0.95 $)]{
	 \includegraphics[scale=0.50]{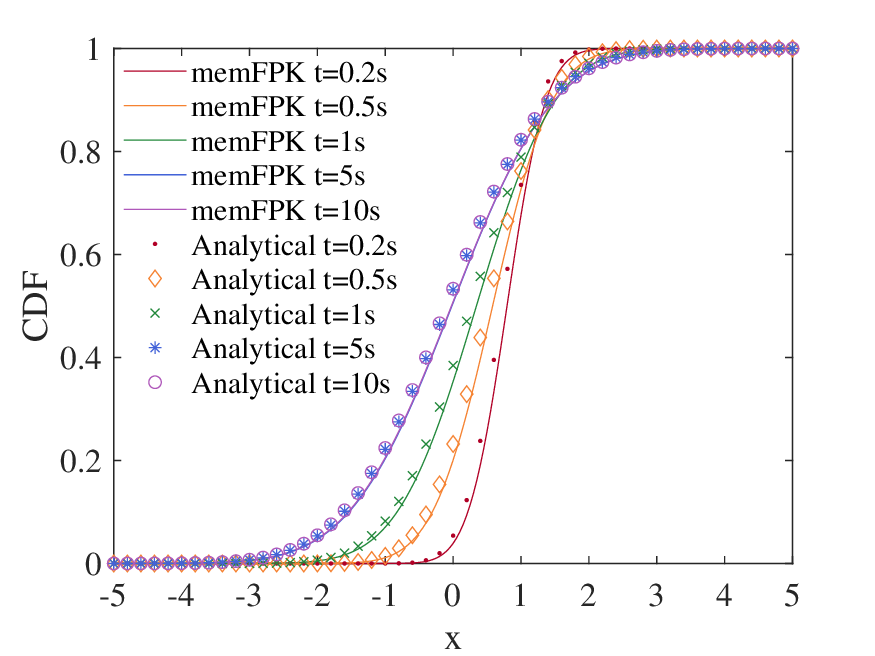}}
\subfloat[logarithmic scale ($ H=0.95 $)]{
	 \includegraphics[scale=0.50]{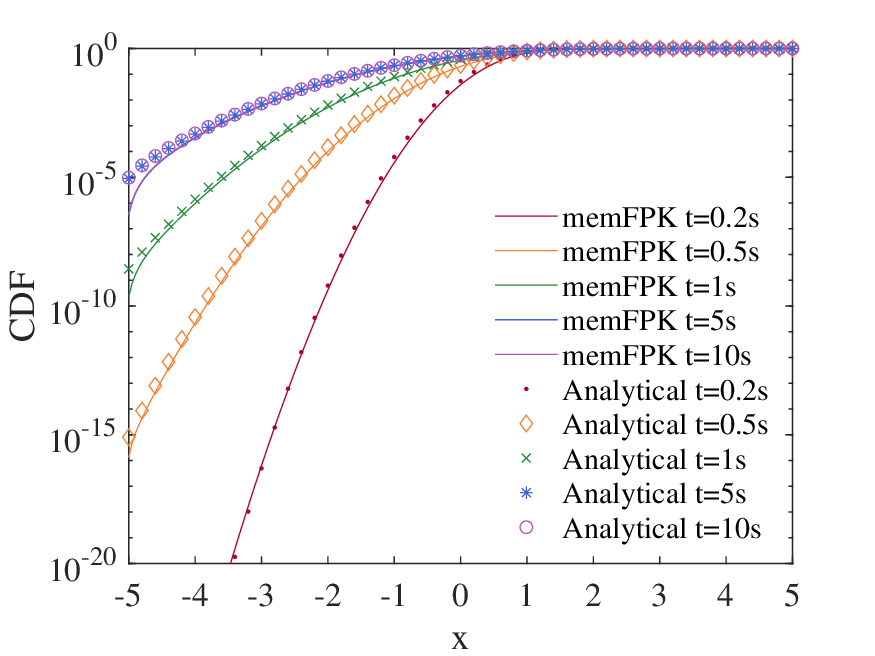}}
\caption{Comparison of the transient CDFs via memFPK equation.}
\label{fig3-1-2}
\end{figure}

In Fig. \ref{fig3-1-1}, the transient PDFs derived from the solutions of the memFPK equation and the analytical solutions are presented for multiple distinct time instants. Concurrently, Fig. \ref{fig3-1-2} showcases the corresponding cumulative distribution functions (CDFs) related to those in Figure \ref{fig3-1-1}. As can be discerned from both Fig. \ref{fig3-1-1} and Fig. \ref{fig3-1-2}, the PDFs and CDFs obtained via the memFPK equation at various time points $t$ are in excellent agreement with the analytical results. This indicates a high level of accuracy, even in the tail regions of the distributions.

\begin{figure}[t!]
\centering
\subfloat[$ H=0.65 $]{
    \includegraphics[scale=0.55]{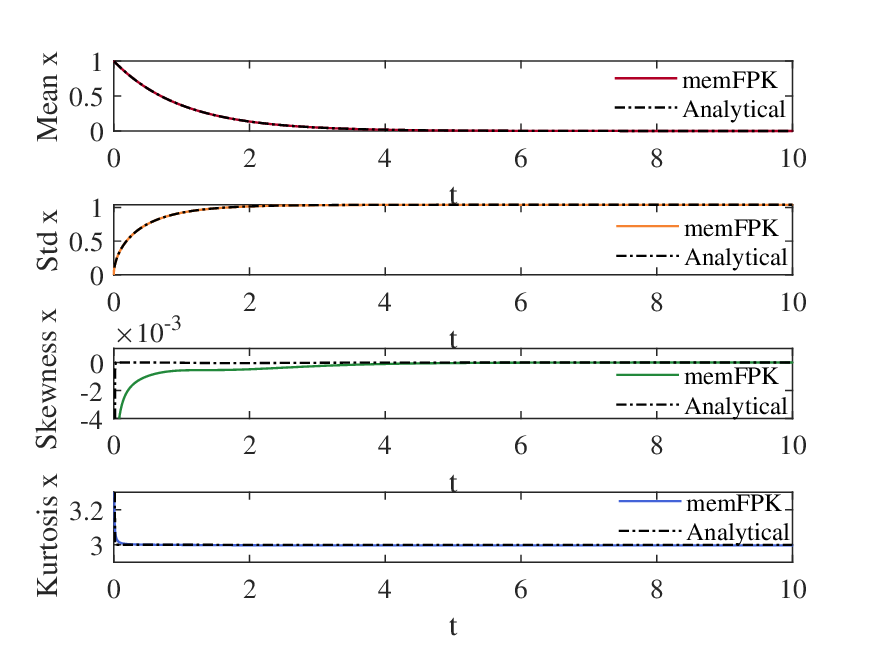}}
\subfloat[$ H=0.95 $]{
	\includegraphics[scale=0.55]{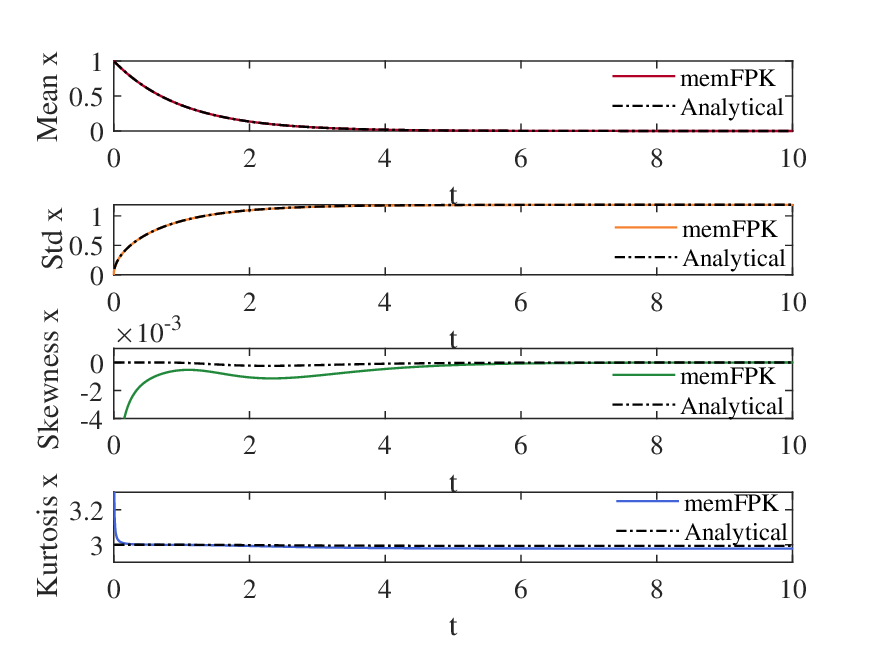}}
\caption{Comparison of the mean, standard deviation, skewness and kurtosis of $ X_t $ via memFPK equation.}
\label{fig3-1-4}
\end{figure}

Fig. \ref{fig3-1-4} depicts the time-varying solutions of the mean, standard deviation, skewness, and kurtosis of the response $X$ obtained from the memFPK equation for two different values of the Hurst parameter $H$. These results are also compared with the corresponding analytical solutions: for $H=0.65$ in Fig. \ref{fig3-1-4}(a) and for $H=0.95$ in Fig. \ref{fig3-1-4}(b). Evidently, the memFPK equation accurately replicates the first four order moments of the response.

\subsection{Duffing oscillator}
Let us consider the over damped bistable duffing oscillator system, governing by Eq. (\ref{mSDS}) as the form		
\begin{equation}\label{example-2}
\dot{X_t}= \alpha X_t+\beta X_t^3+\sigma_W\xi_t+\sigma_B\xi_t^H,
\end{equation}
where $ \alpha, \beta $ are parameters controlling the linear and nonlinear restoring coefficients, the excitations $ \xi_t $ and $ \xi^H_t $ are identified as the unit GWN and FGN, respectively, $ \sigma_W^2 $ and $ \sigma_B^2 $ denote the intensities of the corresponding noise terms. The initial value is taken as $ X_0=x_0 $, a constant. 

Then, refer to Theorem \ref{them-1} (III), the response PDF of Eq. (\ref{example-2}) satisfies the following memFPK equation (\ref{FPKnonlinear}) with the following memory-dependent drift and diffusion coefficients are, 
\begin{align*}
a^{{\rm (mem)}}(x,t)&=\alpha x+\beta x^3,\cr
b^{{\rm (mem)}}(x,t)&=\frac{1}{2}\sigma_W^2+\sigma_B^2\mathbb{E}\Big[\int_{0}^{t}\phi(t,s)\exp\big\{\int_{s}^{t}(\alpha+3\beta X_u^2)\mathrm{d}u\big\}\mathrm{d}s\mid X_t=x\Big],
\end{align*}
respectively. The the memory kernel dependence term
$$ \mathbb{E}\Big[\int_{0}^{t}\phi(t,s)\exp\big\{\int_{s}^{t}(\alpha+3\beta X_u^2)\mathrm{d}u\big\}\mathrm{d}s\mid X_t=x\Big], $$ 
can be approximated by Eq. (\ref{eq3-3}). Then the memFPK equation can be numerically solved via the FD method. Since there is no analytical solution available for Eq. (\ref{example-2}), MCS is carried out. The MCS solution, serving as a reference solution, is utilized to evaluate the effectiveness of the solution of the memFPK equation. The sample paths in MCS are obtained by a 2-stage Runge-Kutta method \cite{hong2021optimal}.

The response interval is $ [-2.5,2.5] $ with the space step $ \Delta x=0.05 $, and the time step $ \Delta t $ takes the value 0.0001 s. {And the parameters of the Duffing oscillator are given in table \ref{table-2}.} 
\begin{table}[h!]
		\centering
		{\caption{Values of parameters of Duffing oscillator for example 5.2}
		\begin{tabular}{ccccccc}
			\toprule
			parameters & $\alpha$ & $\beta$ & $\sigma_W$ & $\sigma_B$ & $x_0$ & $H$ \\ 
			\midrule
			values & 1 & -1 & 0.8 & 0.6 & 0 & 0.6  \\ 
			\bottomrule
		\end{tabular}
		\label{table-2}}
\end{table}

\begin{figure}[t!]
\centering
\subfloat[PDF]{
	\includegraphics[scale=0.50]{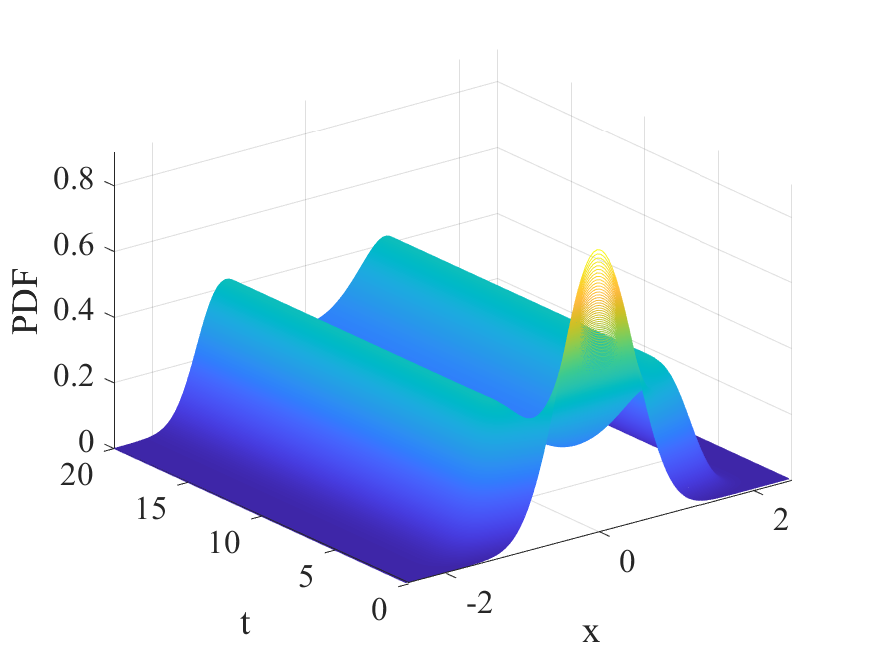}}
\subfloat[MCS PDF]{
	\includegraphics[scale=0.50]{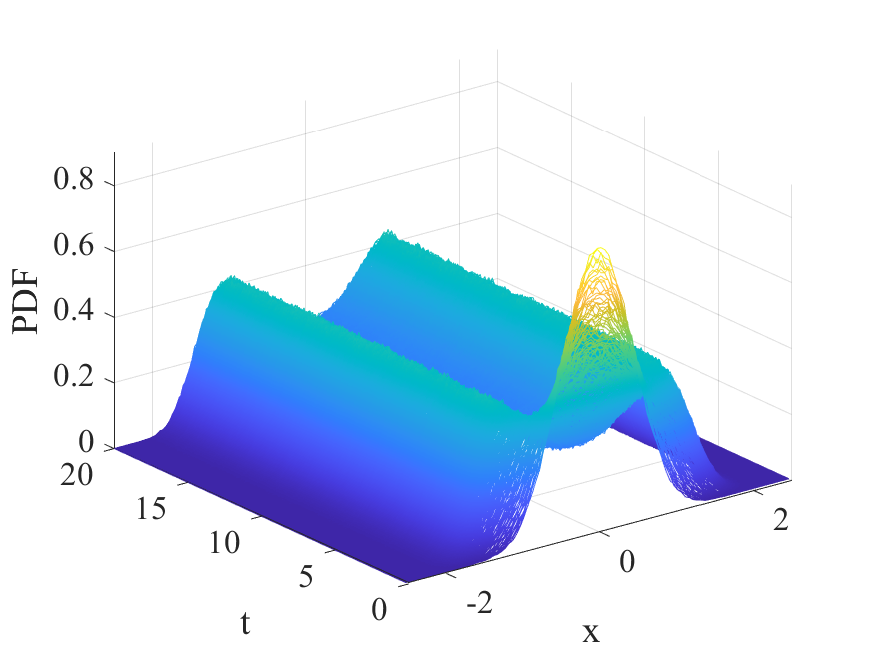}}
\\
\subfloat[Error surface]{
	\includegraphics[scale=0.50]{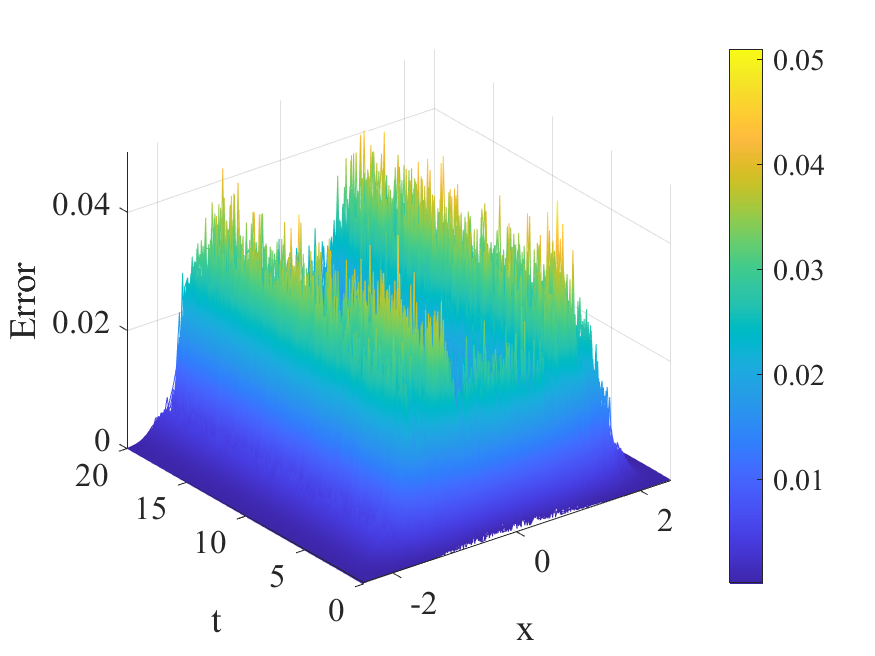}}
\subfloat[Error value]{
	\includegraphics[scale=0.50]{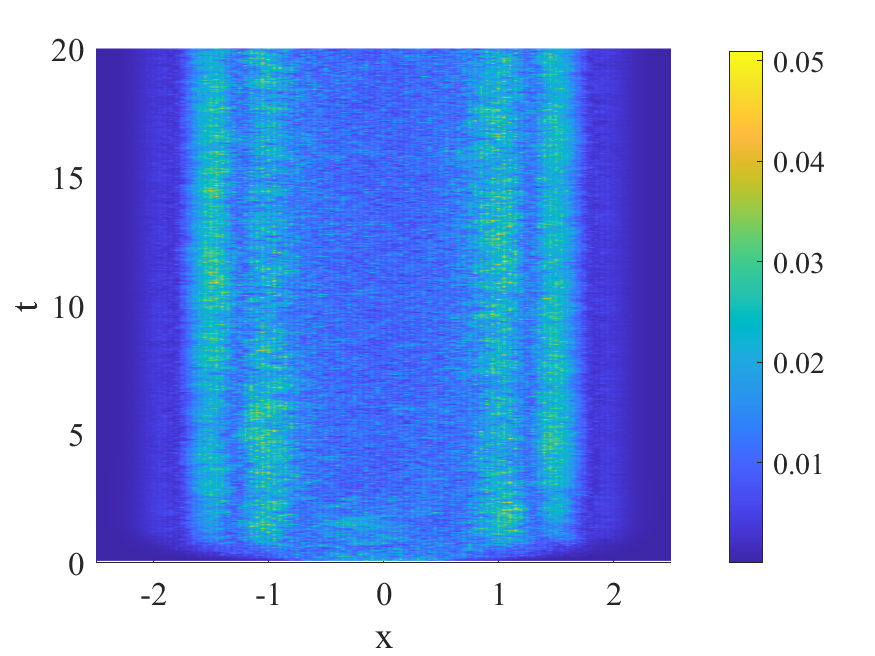}}
\caption{Comparison of the PDF surface via memFPK equation and MCS ($ 10^6 $ samples).}
\label{fig3-2-1}
\end{figure}

Fig. \ref{fig3-2-1} compares the PDF surfaces of $X$ in Eq. (\ref{example-2}). Panel (a) of Fig. \ref{fig3-2-1} shows the results obtained from the memFPK equation, while panel (b) presents the MCS results based on $10^6$ samples. This particular scenario models a strongly nonlinear bistable system under intense random excitation. Panels (c) and (d) of Fig. \ref{fig3-2-1} depict the absolute error between the solutions of the memFPK equation and the MCS results. As evident from Fig. \ref{fig3-2-1}, for the time-varying states of the Duffing oscillator, the solutions of the memFPK equation closely match the MCS results.

\begin{figure}[t!]
\centering
\subfloat[linear scale]{
	\includegraphics[scale=0.50]{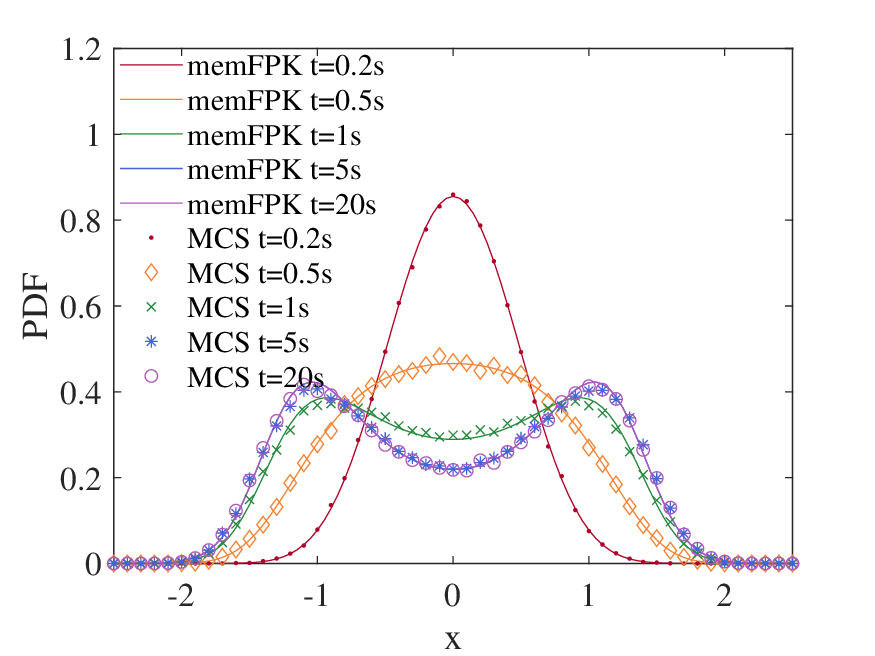}}
\subfloat[logarithmic scale]{
	\includegraphics[scale=0.50]{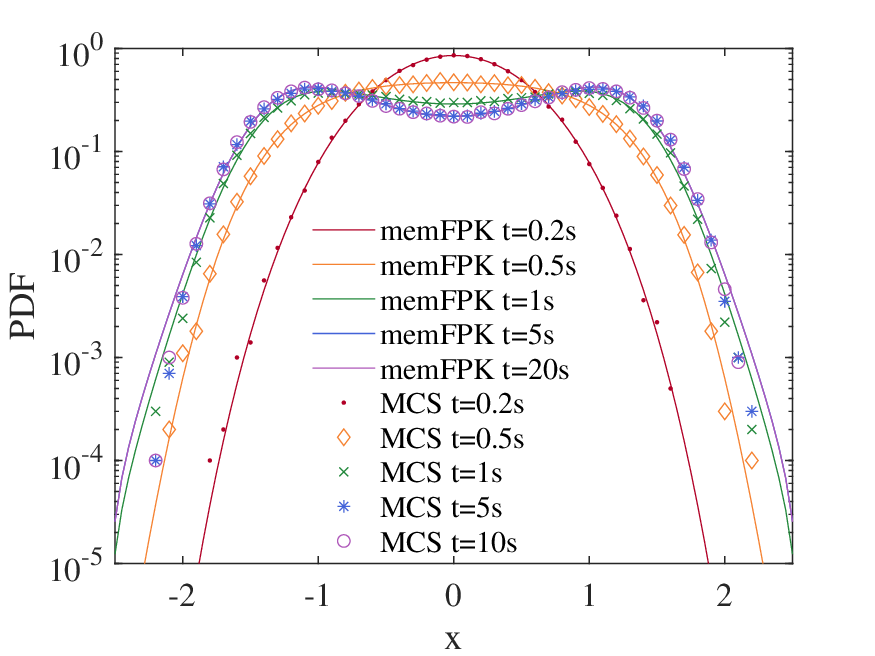}}
\caption{Comparison of the transient PDFs via memFPK equation and MCS  ($ 10^6 $ samples).}
\label{fig3-2-2}
\end{figure}

\begin{figure}[t!]
\centering
\subfloat[linear scale]{
	 \includegraphics[scale=0.50]{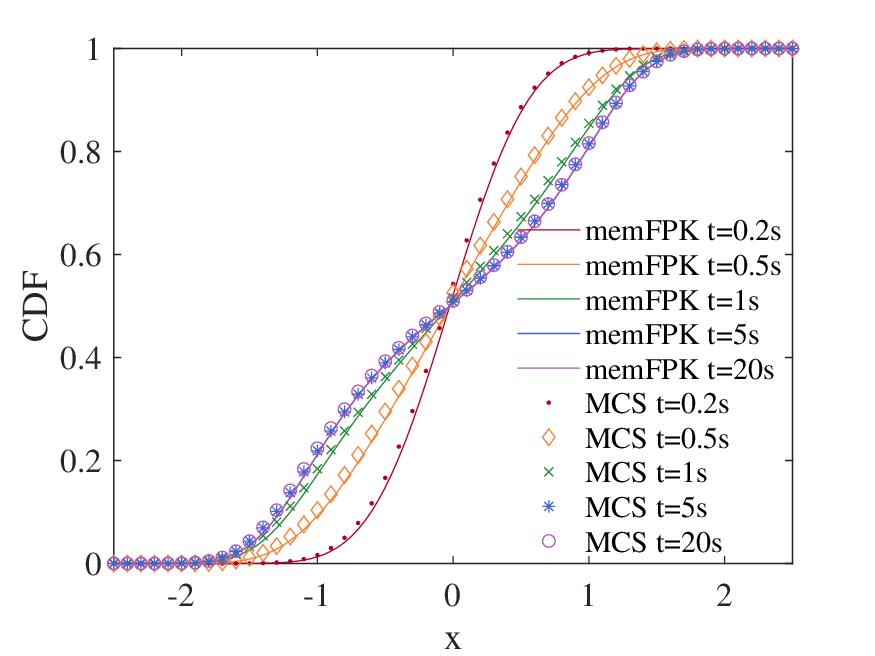}}
\subfloat[logarithmic scale]{
	 \includegraphics[scale=0.50]{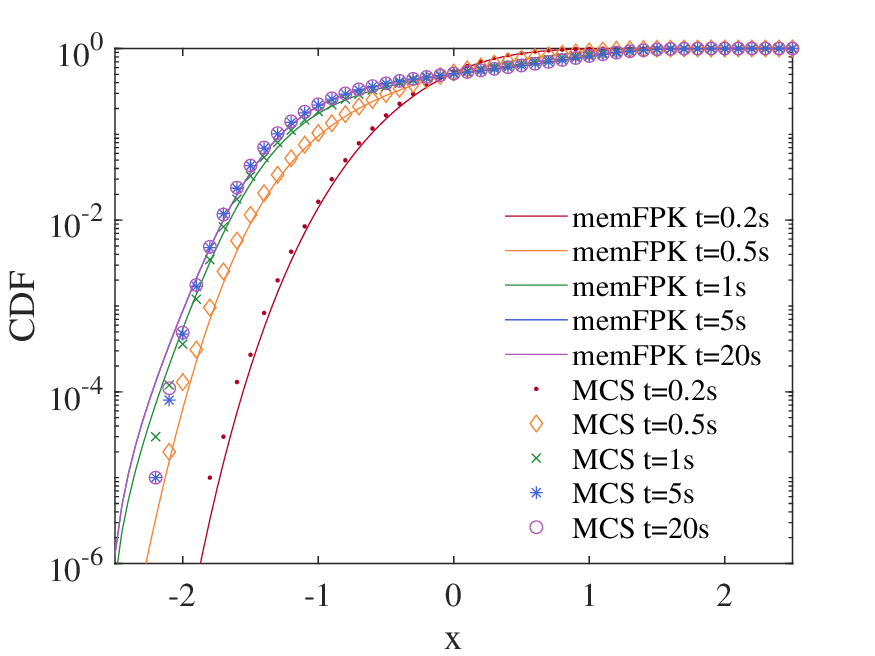}}
\caption{Comparison of the transient CDFs via memFPK equation and MCS  ($ 10^6 $ samples).}
\label{fig3-2-3}
\end{figure}

In addition, Fig. \ref{fig3-2-2} illustrates the PDF curves obtained from the memFPK equation at various time instants, juxtaposed with the corresponding MCS PDFs. The associated CDFs from Fig. \ref{fig3-2-2} are further compared with the MCS results in Fig. \ref{fig3-2-3}. Despite the approximation in Eq. (\ref{eq3-3}) leading to some loss of temporal correlation, the results indicate that the memFPK equation and its approximated form are effective in capturing the time-varying response PDFs induced by non-Markovian excitation. Both the PDF and CDF results exhibit satisfactory accuracy.

\begin{figure}[t!]
\centering
\subfloat[Mean]{
   	\includegraphics[scale=0.50]{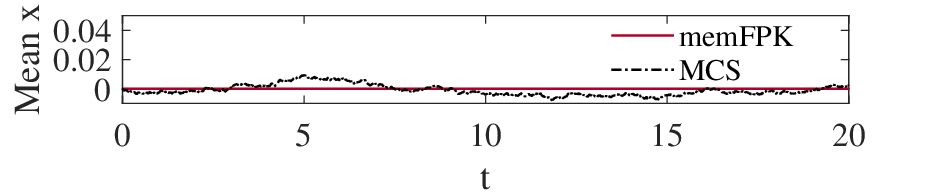}}
\subfloat[Standard deviation]{
   	\includegraphics[scale=0.50]{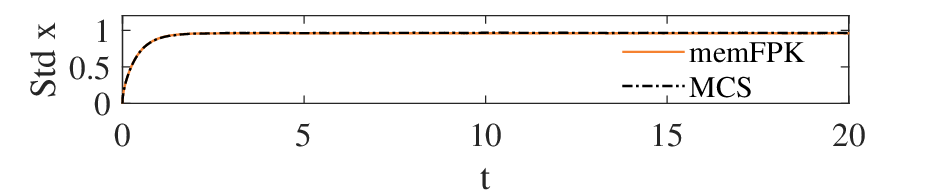}}
\\
\subfloat[Skewness]{
   	\includegraphics[scale=0.50]{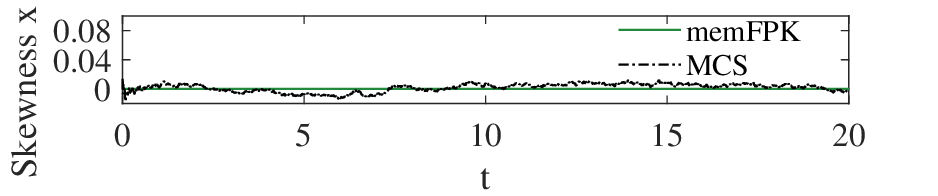}}
\subfloat[Kurtosis]{
   	\includegraphics[scale=0.50]{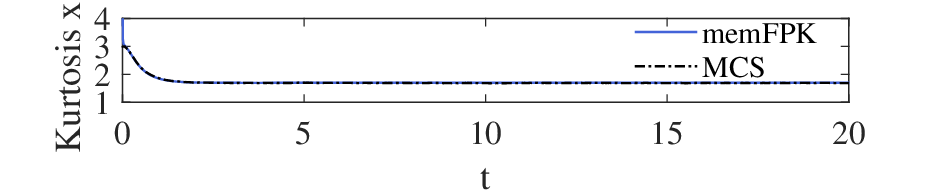}}
\caption{Comparison of the mean, standard deviation, skewness and kurtosis of $ X_t $ via the memFPK equation and the MCS ($ 10^6 $ samples).}
\label{fig3-2-4}
\end{figure}

Moreover, Fig. \ref{fig3-2-4} shows the time-varying solutions for the first four statistical moments (mean, variance, skewness, and kurtosis) of the response 
$X_t$. These moments are computed using the memFPK equation and compared with the MCS solutions. The excellent agreement between the two methods further validates the accuracy and reliability of the memFPK equation in characterizing the dynamic behavior of the system.

\subsection{Verhulst model}
Let us consider the Verhulst model with random environmental disturbances, which describes population dynamics, can be governed by
\begin{equation}\label{example-3}
\dot{X_t}= \alpha X_t-\beta X_t^2+\sigma_WX_t \xi_t+\sigma_BX_t\xi_t^H,
\end{equation}
where $ \alpha $ denotes the intrinsic growth rate in the absence of constraints, $ \beta $ represents the coefficient of intraspecific competition. $ \xi_t $ and $ \xi^H_t $ correspond to the unit GWN and FGN, respectively, while $ \sigma_W^2 $ and $ \sigma_B^2 $ indicate the intensities of the associated noise terms. The initial value $ X_0 $ is a random variable that follows a normal distribution with mean $ \mu_0 $ and variance $ \sigma_0^2 $.

Then, refer to Theorem \ref{them-1} (I), the response PDF of Eq. (\ref{example-3}) satisfies the following memFPK equation (\ref{FPKnonlinear}) with the following memory-dependent drift and diffusion coefficients are, 
\begin{align*}
a^{{\rm (mem)}}(x,t)&=\alpha x-\beta x^2+\frac{1}{2}\sigma_W^2x+\sigma_B^2x\mathbb{E}\Big[\int_{0}^{t}\phi(t,s)\exp\big\{\int_{s}^{t}-2\beta X_u\mathrm{d}u\big\}\mathrm{d}s\mid X_t=x\Big],\cr
b^{{\rm (mem)}}(x,t)&=\frac{1}{2}\sigma_W^2x^2+\sigma_B^2x^2\mathbb{E}\Big[\int_{0}^{t}\phi(t,s)\exp\big\{\int_{s}^{t}-2\beta X_u\mathrm{d}u\big\}\mathrm{d}s\mid X_t=x\Big],
\end{align*}
respectively. Similarly, the memory kernel dependence term $$ \mathbb{E}\Big[\int_{0}^{t}\phi(t,s)\exp\Bigl\{\int_{s}^{t}-2\beta X_u\mathrm{d}u\Bigr\}\mathrm{d}s\mid X_t=x\Big], $$ can be approximated using Eq. (\ref{eq3-3}).
		
The response interval is $ [0,7] $ with the space step $ \Delta x=0.1 $, and the time step $ \Delta t $ takes the value 0.0001 s. {The parameters employed in the Verhulst model are given in table \ref{table-3}.} 
\begin{table}[h!]
	\centering
	{\caption{Values of parameters of Verhulst model for example 5.3}
		\begin{tabular}{cccccccc}
			\toprule
			parameters & $\alpha$ & $\beta$ & $\sigma_W$ & $\sigma_B$ & $\mu_0$ & $\sigma_0^2$ & $H$ \\ 
			\midrule
			values & 4 & 1 & 0.1 & 0.3 & 1 & 0.05 & 0.8  \\ 
			\bottomrule
		\end{tabular}
		\label{table-3}}
\end{table}

\begin{figure}[t!]
\centering
\subfloat[PDF]{
	\includegraphics[scale=0.50]{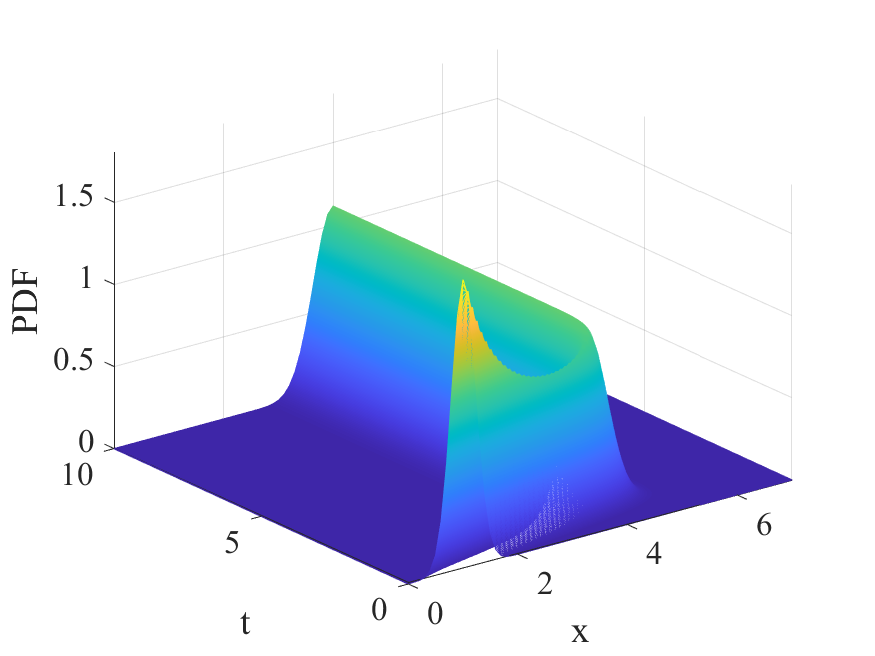}}
\subfloat[MCS PDF]{
	\includegraphics[scale=0.50]{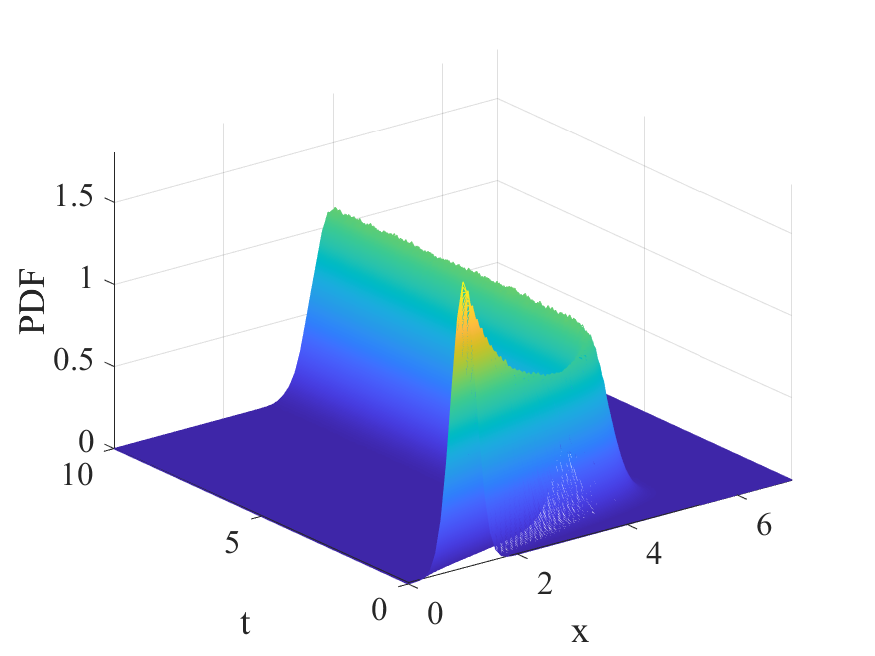}}
\\
\subfloat[Error surface]{
	\includegraphics[scale=0.50]{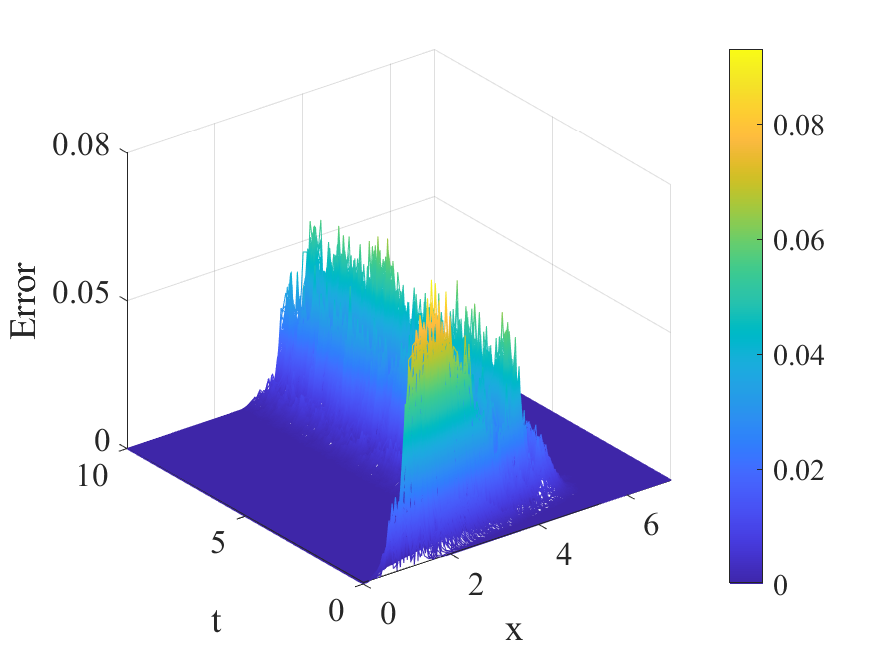}}
\subfloat[Error value]{
	\includegraphics[scale=0.50]{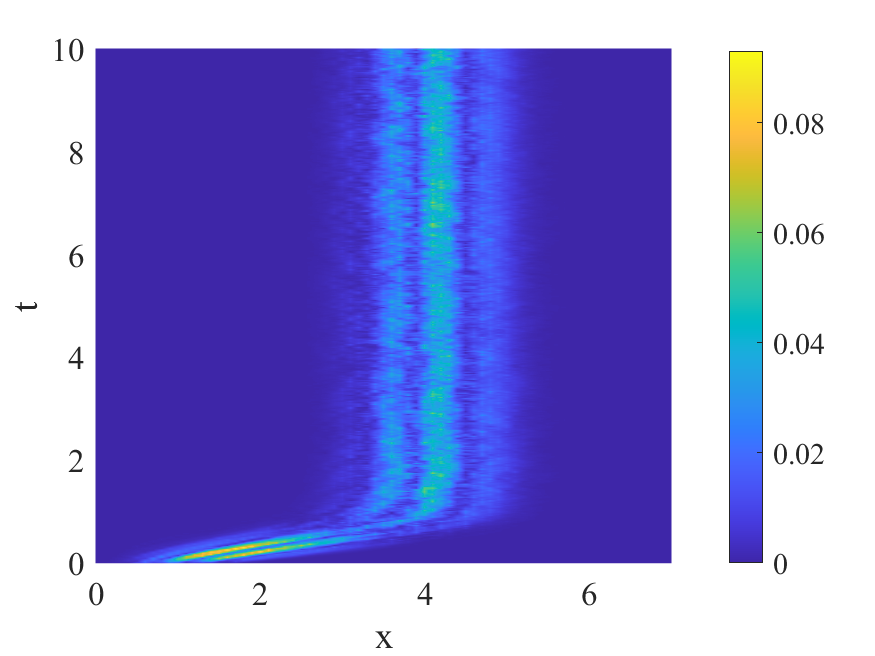}}
\caption{Comparison of the PDF surface via memFPK equation and MCS ($ 10^6 $ samples).}
\label{fig3-3-1}
\end{figure}

\begin{figure}[t!]
\centering
\subfloat[linear scale]{
	\includegraphics[scale=0.50]{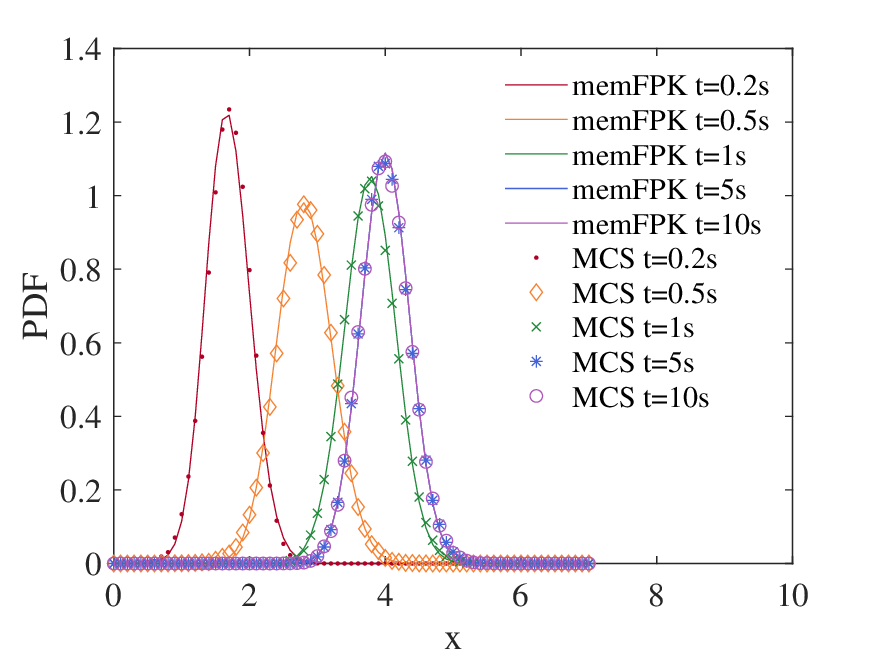}}
\subfloat[logarithmic scale]{
	\includegraphics[scale=0.50]{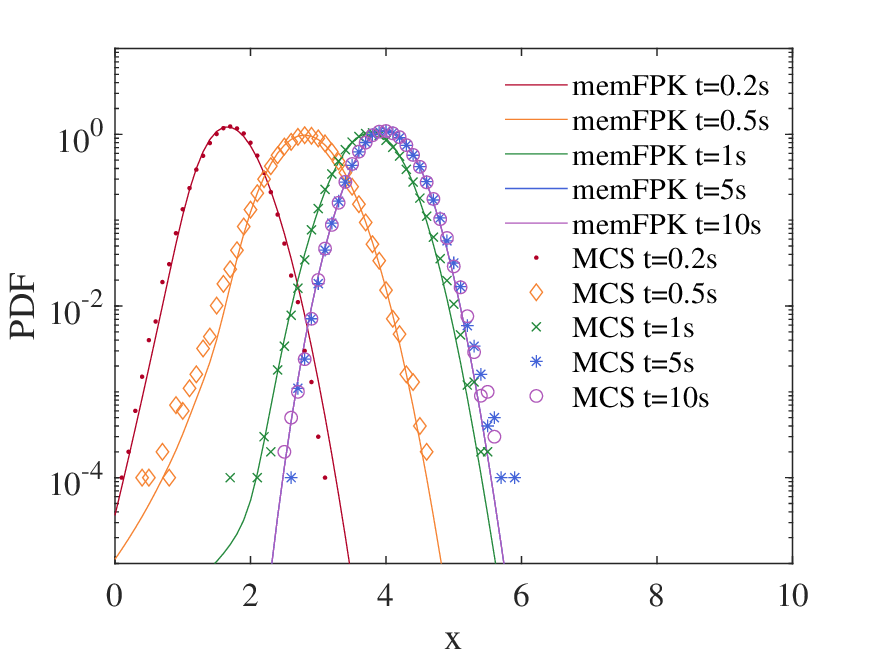}}
\caption{Comparison of the transient PDFs via memFPK equation and MCS ($ 10^6 $ samples).}
\label{fig3-3-2}
\end{figure}

\begin{figure}[t!]
\centering
\subfloat[linear scale]{
    \includegraphics[scale=0.50]{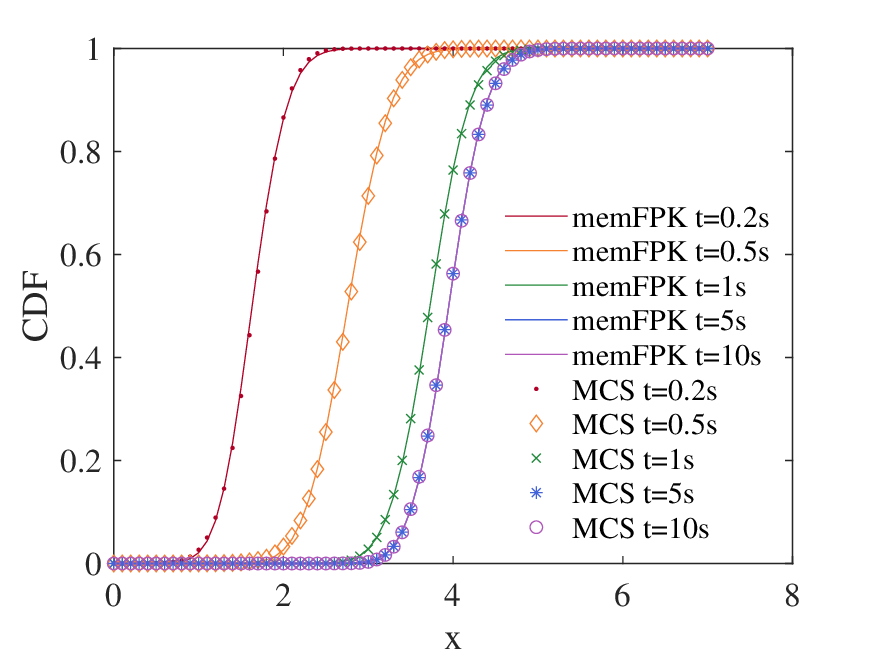}}
\subfloat[logarithmic scale]{
    \includegraphics[scale=0.50]{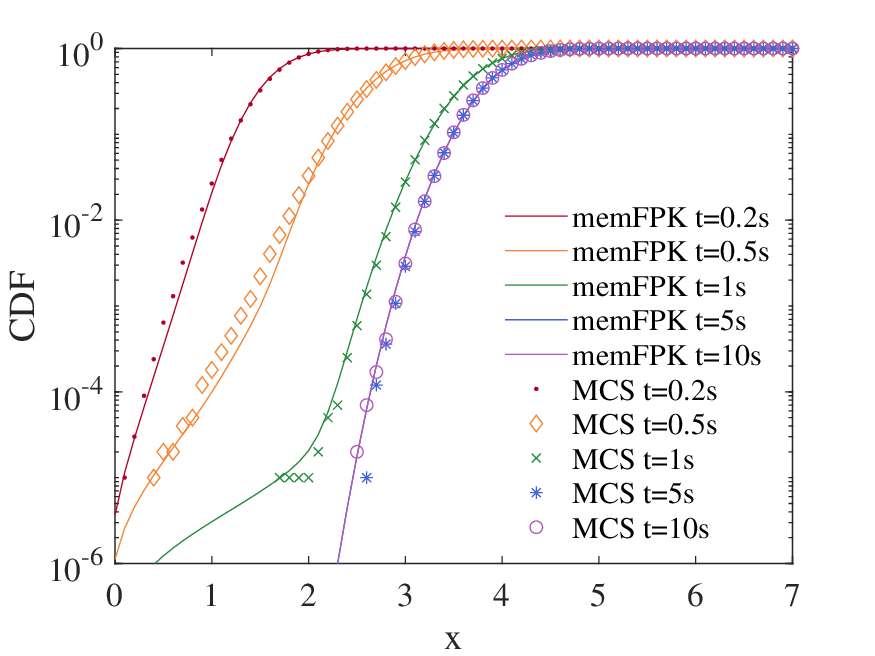}}
\caption{Comparison of the transient CDFs via memFPK equation and MCS ($ 10^6 $ samples).}
\label{fig3-3-3}
\end{figure}

Fig. \ref{fig3-3-1} displays the solutions of the memFPK equation and MCS results for the PDF surface of 
$X_t$ in Eq. (\ref{example-3}), along with the corresponding absolute error. Meanwhile, Figs. \ref{fig3-3-2} and \ref{fig3-3-3} illustrate the PDFs and CDFs, respectively, obtained from the memFPK equation and MCS solutions at several representative time instants $t$. This case models a strongly nonlinear system under random excitation. As evident from Fig. \ref{fig3-3-1}, the solutions of the memFPK equation show excellent agreement with the MCS results throughout the considered scenario.

\begin{figure}[t!]
\centering
\subfloat[Mean]{
    \includegraphics[scale=0.50]{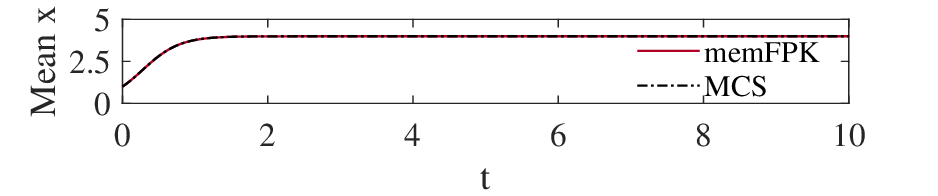}}
\subfloat[Standard deviation]{
    \includegraphics[scale=0.50]{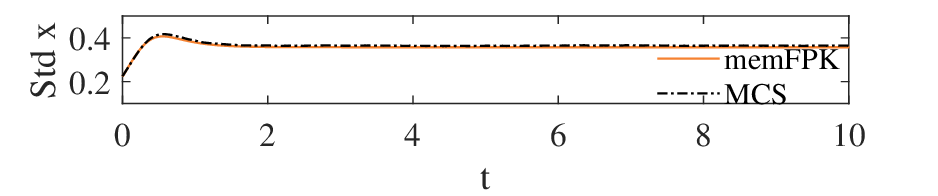}}
    \\
\subfloat[Skewness]{
    \includegraphics[scale=0.50]{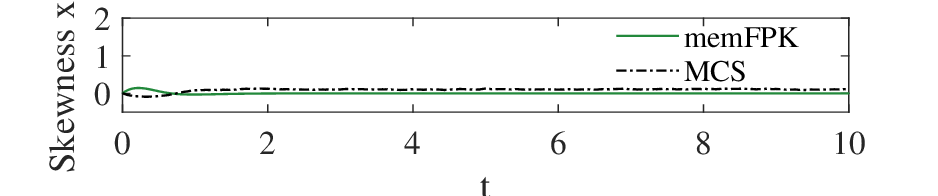}}
\subfloat[Kurtosis]{
    \includegraphics[scale=0.50]{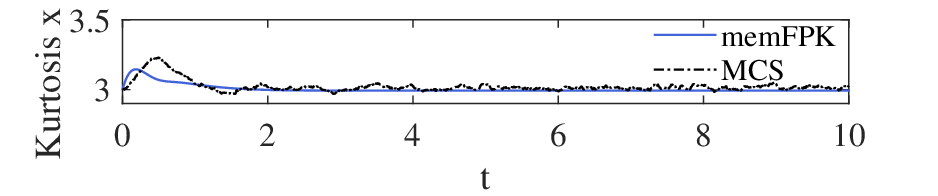}}
\caption{Comparison of the mean, standard deviation, skewness and kurtosis of $ X_t $ via memFPK equation and MCS ($ 10^6 $ samples).}
\label{fig3-3-4}
\end{figure}

Moreover, Fig. \ref{fig3-3-4} presents the time-varying solutions for the first four-order moments of the response $X_t$, namely the mean, variance, skewness, and kurtosis. These moments are calculated using the memFPK equation and compared with the MCS solutions. The results demonstrate a high level of consistency between the two approaches, further validating the effectiveness of the memFPK equation in characterizing the dynamic behavior of the system under stochastic excitation.
	
\subsection{Two-DOF quasi-non-integrable Hamiltonian system}
Consider the following two-DOF quasi-non-integrable Hamiltonian system
\begin{align}\label{example-4-0}
&\ddot{X_1}+\gamma\dot{X_1}+\frac{\partial U(X_1,X_2)}{\partial X_1}=\sqrt{2D_1}\xi_t^{1,H},\cr
&\ddot{X_2}+\gamma\dot{X_2}+\frac{\partial U(X_1,X_2)}{\partial X_2}=\sqrt{2D_2}\xi_t^{2,H},
\end{align}
where $ \gamma $ is a constant coefficient of linear damping, $ \xi_t^{1,H} $ and $ \xi_t^{2,H} $ are two independent unit FGNs, $ 2D_1 $, $ 2D_2 $ are the intensities of excitations.
		
Let $ Q_1 = X_1, Q_2 = X_2, P_1 = \dot{X_1}, P_2 =  \dot{X_2} $, the Hamiltonian for system (\ref{example-4-0}) is
\begin{align}\label{example-4-1}
&\mathcal{H} = P_1^2/2 + P_2^2/2 + U(Q_1, Q_2),\cr
&U(Q_1, Q_2)=\omega_1^2Q_1^2/2+\omega_2^2Q_2^2/2+\lambda(\omega_1^2Q_1^2/2+\omega_2^2Q_2^2/2)^2/4,
\end{align}
where $ \omega_1, \omega_2 $ are the angle frequencies of linearized system; $ \lambda $ indicates the intensity of nonlinearity.
		
Deng et al. \cite{deng2016stochastic} use the stochastic averaging method for quasi-non-integrable Hamiltonian systems under FGN excitation, resulting in the following averaged one-dimensional SDE
\begin{align}\label{example-4}
\mathrm{d}\mathcal{H} = m(\mathcal{H})\mathrm{d}t+\sigma(\mathcal{H})\mathrm{d}^{\circ}B^H_t,
\end{align}
where the initial value is taken as $ \mathcal{H}_0=h_0 $, a constant, and the drift and diffusion coefficients are as follows
\begin{align}\label{example-4-2}
&m(\mathcal{H})=-2\gamma\big(\mathcal{H}-\frac{1}{4}R+\frac{\lambda}{12}R^2\big),\cr
&\sigma(\mathcal{H})=\Big[2(D_1+D_2)\big(\mathcal{H}-\frac{1}{4}R+\frac{\lambda}{12}R^2\big)\Big]^{1/2},\cr
&R=\frac{1}{\lambda}(\sqrt{1+4\lambda\mathcal{H}}-1).
\end{align}
		
Then, refer to Theorem \ref{them-1} (I), the response PDF of SDE (\ref{example-4-1}) satisfies the following memFPK equation (\ref{FPKnonlinear}) with the following memory-dependent drift and diffusion coefficients are, 
\begin{align*}
a^{{\rm (mem)}}(h,t)&=\frac{6\gamma}{\lambda}(-1-8\lambda h+\sqrt{1+4\lambda h})+\frac{(D_1+D_2)(-1+4\sqrt{1+4\lambda x})}{6\sqrt{1+4\lambda x}}\cr
&\quad\times\mathbb{E}\Big[\int_{0}^{t}\phi(t,s)\exp\big\{\int_{s}^{t}
\frac{\gamma}{6}\Big(-4+\frac{1}{\sqrt{1+4\lambda \mathcal{H}_u}}\Big)\mathrm{d}u\Big\}\mathrm{d}s\mid \mathcal{H}_t=h\big],\cr
b^{{\rm (mem)}}(h,t)&=\frac{(D_1+D_2)(1+8\lambda h-\sqrt{1+4\lambda h})}{6\lambda}\cr
&\quad\times\mathbb{E}\Big[\int_{0}^{t}\phi(t,s)\exp\big\{\int_{s}^{t}\frac{\gamma}{6}\Big(-4+\frac{1}{\sqrt{1+4\lambda \mathcal{H}_u}}\big)\mathrm{d}u\big\}\mathrm{d}s\mid \mathcal{H}_t=h\Big],
\end{align*}
respectively. The memory kernel dependence term $$ \mathbb{E}\Big[\int_{0}^{t}\phi(t,s)\exp\big\{\int_{s}^{t}\frac{\gamma}{6}\Big(-4+\frac{1}{\sqrt{1+4\lambda \mathcal{H}_u}}\big)\mathrm{d}u\big\}\mathrm{d}s\mid \mathcal{H}_t=h\Big], $$ is approximated by Eq. (\ref{eq3-3}).

The response interval is $ [0,10] $ with the space step $ \Delta h=0.1 $, and the time step $ \Delta t $ takes the value 0.0001 s. {The parameters employed in two-DOF quasi-non-integrable Hamiltonian system are given in table \ref{table-4}.} 
\begin{table}[h!]
	\centering
	{\caption{Values of parameters of two-DOF quasi-non-integrable Hamiltonian system for example 5.4}
		\begin{tabular}{ccccccccc}
			\toprule
			parameters & $\omega_1$ & $\omega_2$ & $\lambda$ & $\gamma$ & $D_1$ & $D_2$ & $h_0$ & $H$ \\ 
			\midrule
			values & 1.414 & 2 & 2 & 0.01 & 0.01 & 0.01 & 2 & 0.7 \\ 
			\bottomrule
		\end{tabular}
		\label{table-4}}
\end{table}

\begin{figure}[t!]
\centering
\subfloat[PDF]{
	\includegraphics[scale=0.50]{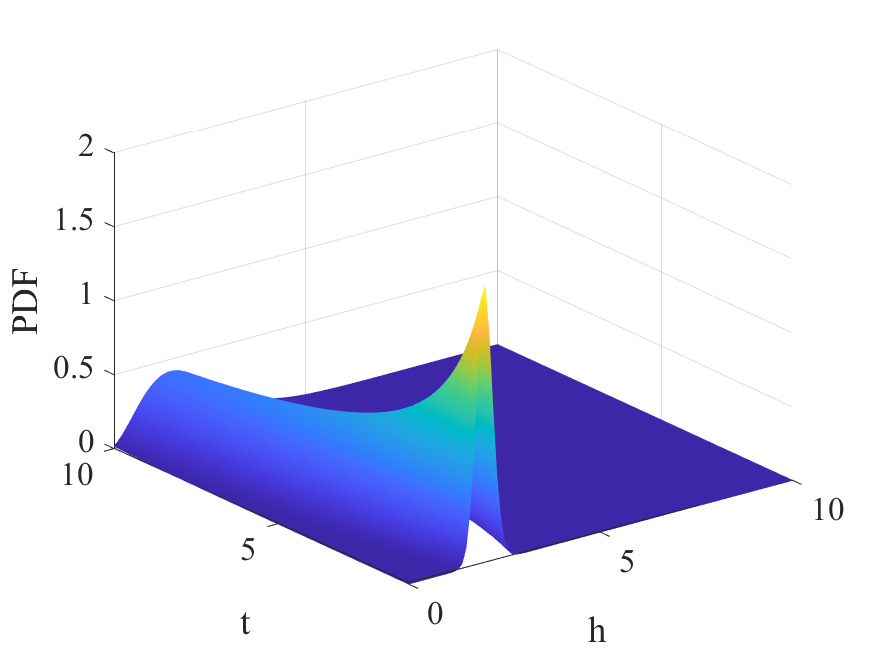}}
\subfloat[MCS PDF]{
	\includegraphics[scale=0.50]{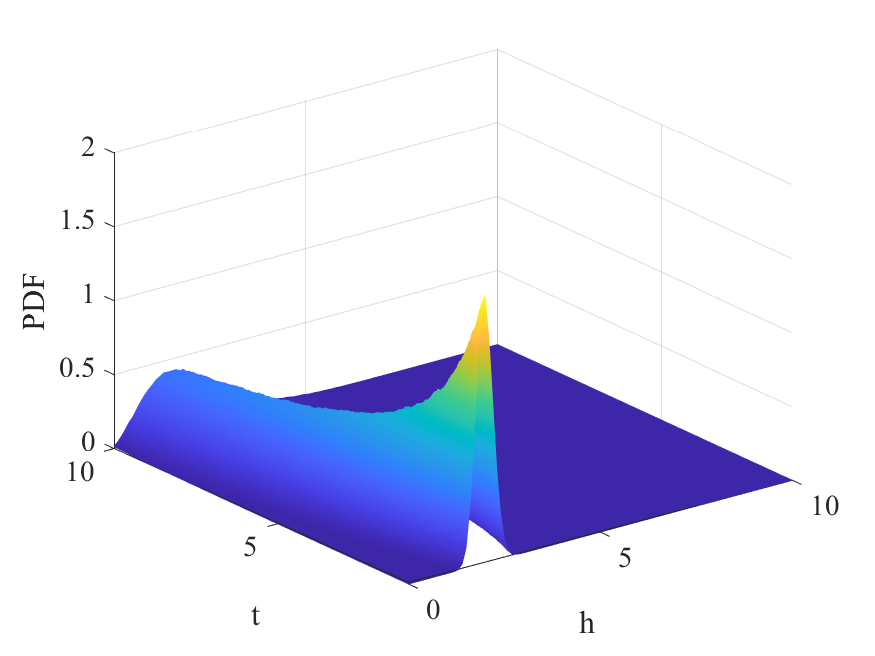}}
\\
\subfloat[Error surface]{
	\includegraphics[scale=0.50]{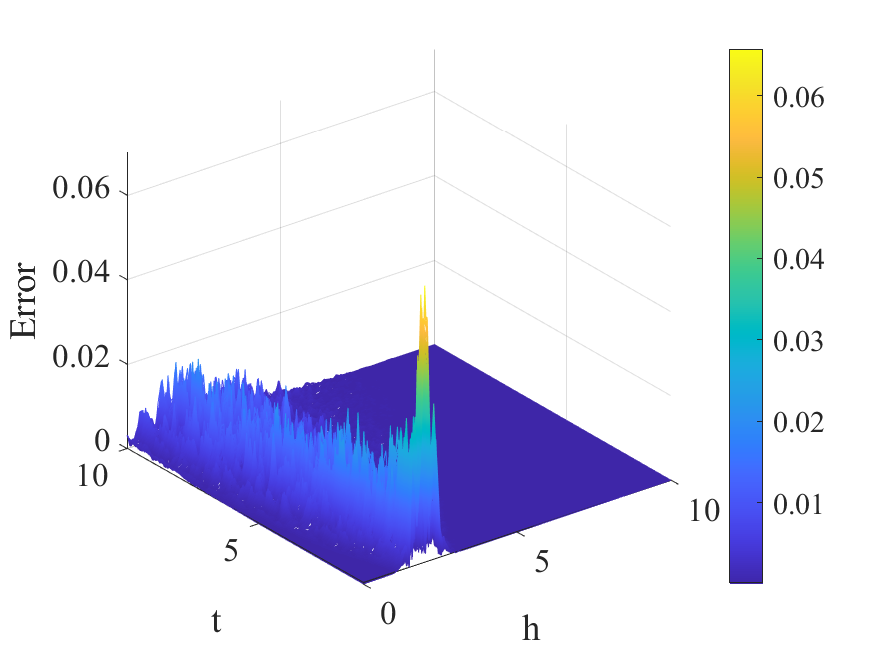}}
\subfloat[Error value]{
	\includegraphics[scale=0.50]{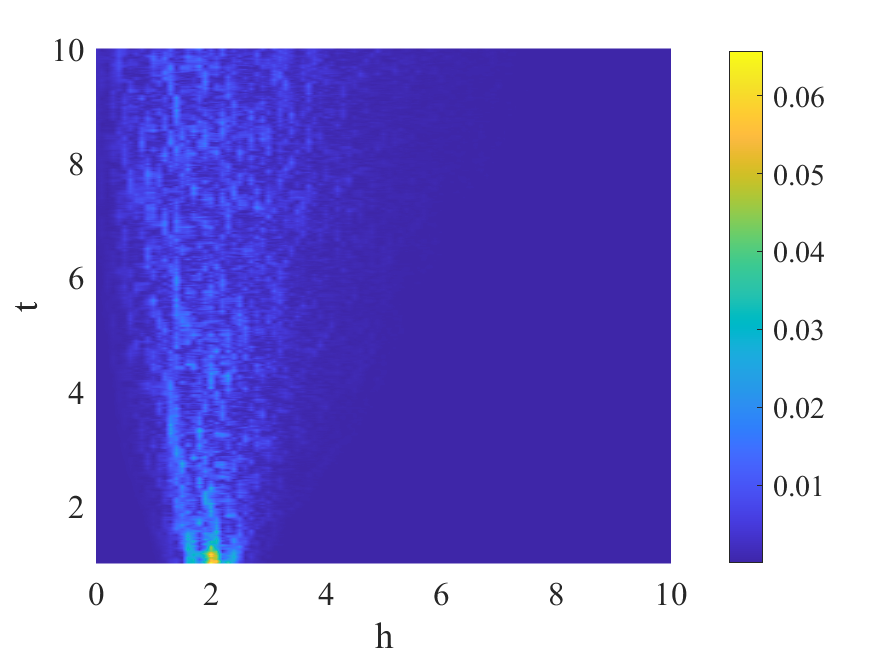}}
\caption{Comparison of the PDF surface via memFPK equation and MCS ($ 10^6 $ samples).}
\label{fig3-4-1}
\end{figure}

\begin{figure}[t!]
\centering
\subfloat[linear scale]{
	\includegraphics[scale=0.50]{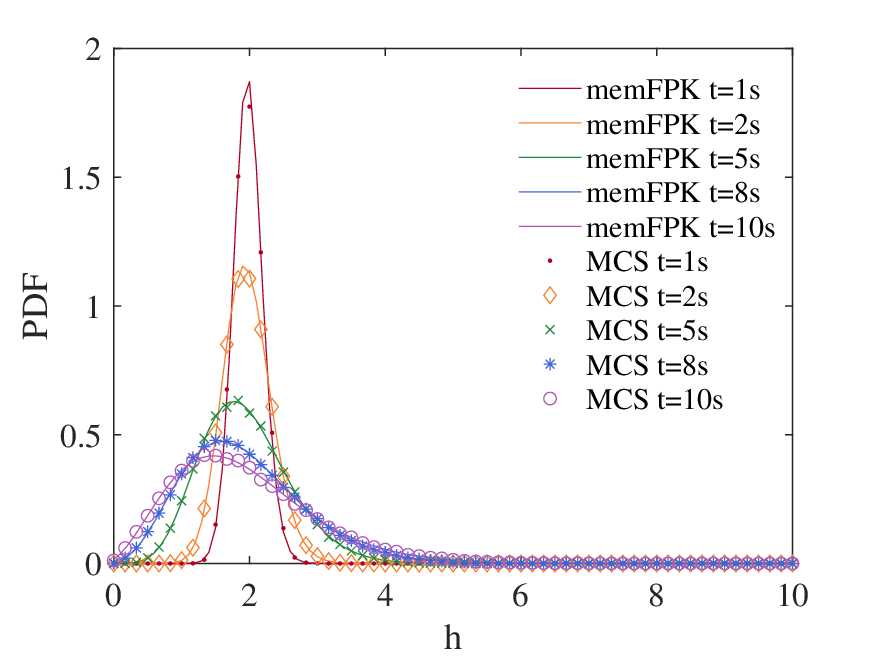}}
\subfloat[logarithmic scale]{
	\includegraphics[scale=0.50]{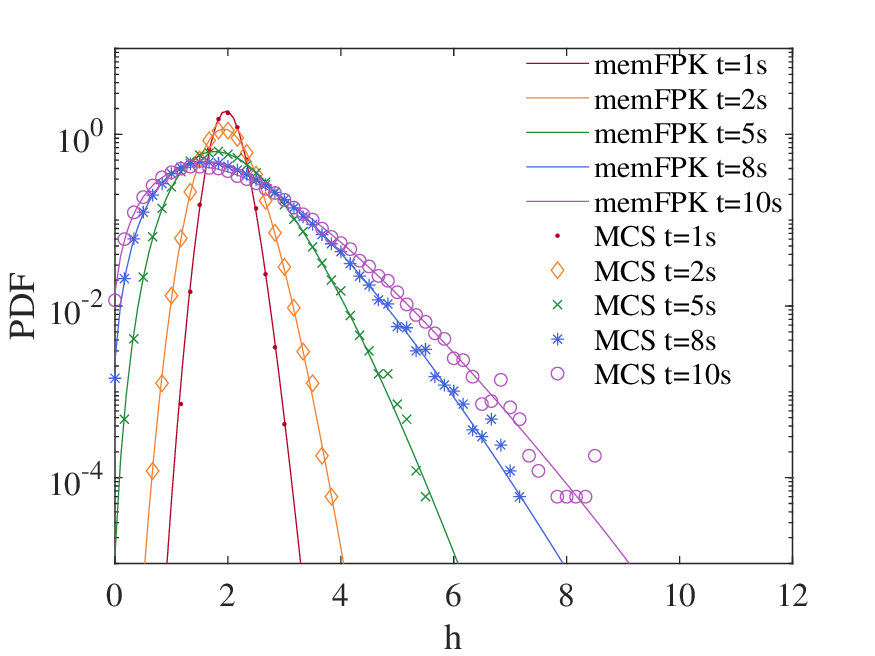}}
\caption{Comparison of the transient PDFs via memFPK equation and MCS ($ 10^6 $ samples).}
\label{fig3-4-2}
\end{figure}

\begin{figure}[t!]
\centering
\subfloat[linear scale]{
	\includegraphics[scale=0.50]{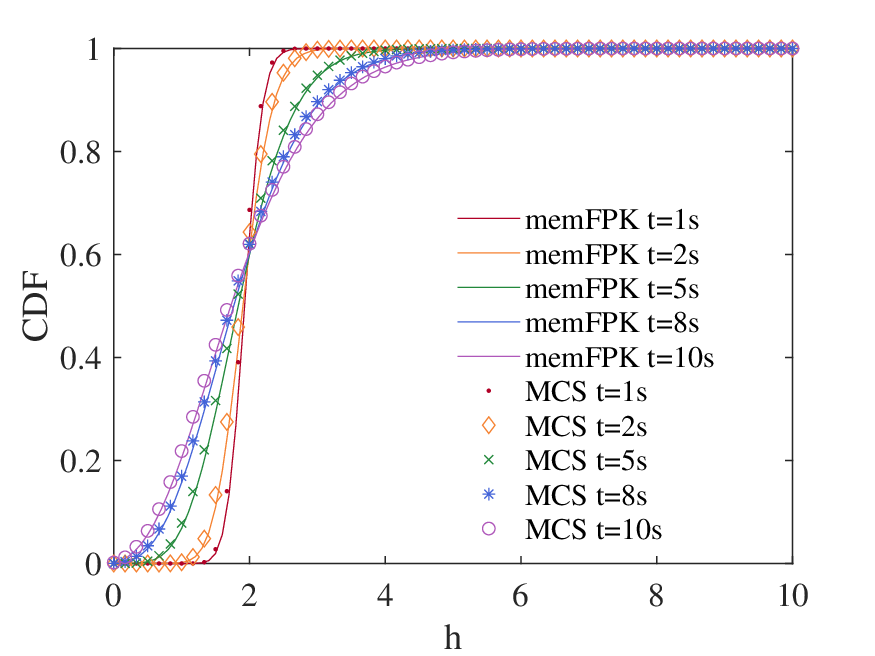}}
\subfloat[logarithmic scale]{
	\includegraphics[scale=0.50]{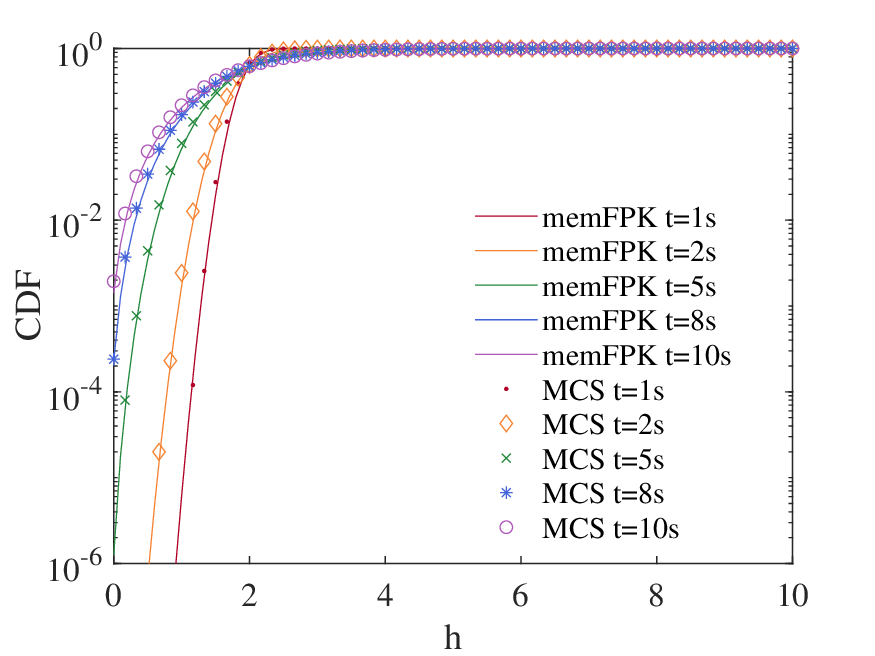}}
\caption{Comparison of the transient CDFs via memFPK equation and MCS ($ 10^6 $ samples).}
\label{fig3-4-3}
\end{figure}

\begin{figure}[t!]
\centering
\subfloat[Mean]{
	\includegraphics[scale=0.50]{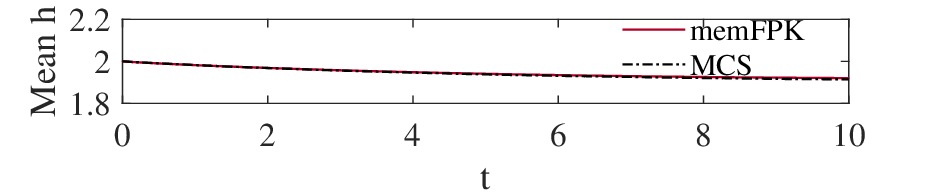}}
\subfloat[Standard deviation]{
	\includegraphics[scale=0.50]{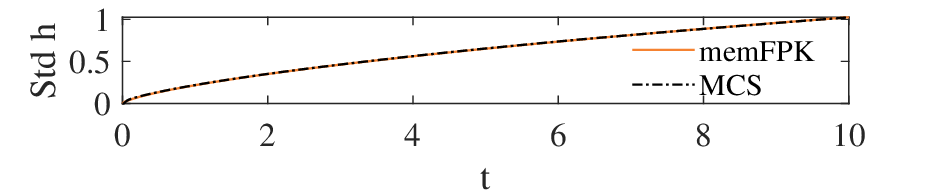}}
\\
\subfloat[Skewness]{
	\includegraphics[scale=0.50]{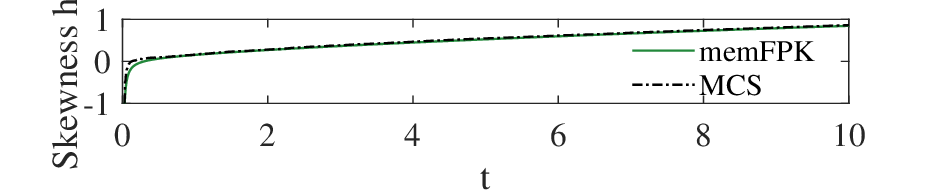}}
\subfloat[Kurtosis]{
	\includegraphics[scale=0.50]{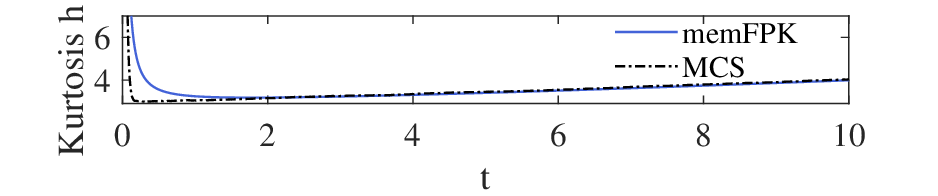}}
\caption{Comparison of the mean, standard deviation, skewness and kurtosis of $ \mathcal{H}_t $ via memFPK equation and MCS ($ 10^6 $ samples).}
\label{fig3-4-4}
\end{figure}

The solutions for PDF and CDF of $ \mathcal{H}_t $ in SDE (\ref{example-4}), obtained from both the memFPK equation and MCS, are presented in Figs. \ref{fig3-4-1}-\ref{fig3-4-3}. Additionally, Fig. \ref{fig3-4-4} compares the time evolution of the first fourth order moments (mean, variance, skewness, and kurtosis) of $ \mathcal{H}_t $. The results demonstrate a good agreement.
				
\section{Discussions and conclusions}\label{sec13}

In this paper, a new memFPK equation for probabilistic responses of nonlinear dynamical systems driven by FGN and GWN is proposed using FWIS integral theory.  Specifically, the historical dependence of the non-Markovian dynamical system's state, which arises from the FGN excitation terms, is accounted for by incorporating the memory kernel dependence term into the coefficients of memFPK equation. Significantly, the developed method can also be regarded as an extension of the findings in Ref. \cite{pei2025non}, allowing for the analysis of non-Markovian dynamical systems with more general drift and diffusion coefficients.
To demonstrate the reliability of the memFPK equation, several numerical examples related to mechanics and dynamics have been investigated. These examples encompass scenarios with both additive and multiplicative noise excitations, as well as linear and nonlinear systems. Comparisons with relevant analytical and MCS solutions have shown that the proposed technique exhibits a high degree of accuracy. Given the pivotal role of PDFs in non-Markovian dynamical analysis, the development of the memFPK equation lays a valuable foundation for future research on nonlinear dynamical systems driven by combined FGN and GWN. {Moreover, the memFPK equation holds potential extensions for broader applications in applied mechanics, particularly in control and reliability analysis of nonlinear systems under long-memory excitations. Additionally, it should be emphasized that owing to the FWIS integration theory, the proposed memFPK framework is not restricted to one-dimensional cases but can be feasible to extended to general multi-dimensional nonlinear systems subjected to both additive and multiplicative FGN excitations, though it may introduce additional mathematical and computational challenges. Addressing these challenges, this line of research is currently in progress.}

\section*{CRediT authorship contribution statement}
\textbf{Lifang Feng:} Conceptualization, Formal analysis, Investigation, Methodology, Software, Visualization, Validation, Writing-original draft, Writing-review \& editing.
\textbf{Bin Pei:} Conceptualization, Funding acquisition, Methodology, Writing-original draft, Writing-review \& editing.
\textbf{Yong Xu:} Conceptualization, Funding acquisition, Project administration, Supervision, Writing-review \& editing.

\section*{Declaration of competing interest}
The authors declare that they have no known competing financial interests or personal relationships that could have appeared to influence the work reported in this paper.

\section*{Acknowledgements}
B. Pei and L. Feng  were supported by National Natural Science Foundation (NSF) of China under Grant No. 12172285, Guangdong Basic and Applied Basic Research Foundation under Grant No. 2024A1515012164, and the Fundamental Research Funds for the Central Universities. Y. Xu was supported by the NSF of China under Grant No. U2441204.		 
		
\appendix
\renewcommand{\thetheorem}{\Alph{section}.\arabic{theorem}}
\section{Some useful Lemmas}\label{app-sec-1}
Let $L^p(0, T), p\geq1$ be the space of all scalar-valued, measurable stochastic processes $(F_t,t\in[0,T])$ such that $ \mathbb{E}[\int_{0}^{T}|F_t|^p\mathrm{d}t]<\infty$. And denote $\mathscr{L}^2_{\phi}(0, T)$ be be the space of measurable, scalar-valued stochastic processes $(F_t,t\in[0,T])$, such that
$ \mathbb{E}[\int_0^T\int_0^T F_s F_t \phi(s,t)\mathrm{d}s\mathrm{d}t\Big]<\infty$ with $\phi(t,s)=H(2H-1)|t-s|^{2H-2}$. By employing similar arguments as in \cite[Theorem 4.3]{ducan2000stochastic} and \cite[Theorem 3.1]{sonmez2023mixed}, the following lemmas can be derived, which are essential for obtaining the fractional It\^o formula for both It\^o-type and Stratonovich-type integrals.
\begin{lemma}\label{itomixed}
Let $ X_t $ be the solution of (\ref{mixSDE}). If $f(X_t)\in L^1(0, T), g(X_t)\in L^2(0, T)$ and $ h(X_t)\in \mathscr{L}^2_{\phi}(0, T)$. Assume that there is an $\alpha>1-H $ such that
\begin{align*}
\mathbb{E}\big[|h(X_u)-h(X_v)|\big]\leq C|u-v|^{2\alpha},
\end{align*}
where $ |u-v| \leq \delta $ for some $ \delta > 0 $ and
\begin{align*}
\lim_{0\leq u,v\leq t,|u-v|\to 0}\mathbb{E}\big[|D^{\phi}_u(h(X_u)-h(X_v))|^2\big]\leq C|u-v|^{2\alpha}.
\end{align*}
			
Then, suppose that $ F(t,x) $ is any function satisfying continuously twice differentiable in $x$ and once in $t$. Moreover, it is assumed that $ h(X_t)D^{\phi}_tX_t\in L^1(0, T)$ and $ h(X_t)F_x(t,X_t)\in \mathscr{L}^2_{\phi}(0, T)$. Then for $ t\in [0,T] $
\begin{align}\label{itofor}
\mathrm{d}F(t,X_t)=&\big\{\frac{\partial F}{\partial t}(t,X_t)+\big(f(X_t)+\frac{1}{2}g(X_t)g'(X_t)+h'(X_t)D^{\phi}_tX_t\big)\frac{\partial F}{\partial x}(t,X_t)\cr
&+\big(\frac{1}{2}g^2(X_t)+h(X_t)D^{\phi}_tX_t\big)\frac{\partial^2 F}{\partial x^2}(t,X_t)\big\}\mathrm{d}t\cr
&+g(X_t)\frac{\partial F}{\partial x}(t,X_t)\mathrm{d}W_t+h(X_t)\frac{\partial F}{\partial x}(t,X_t) \mathrm{d}^{\diamond}B^H_t,\quad a.s.,
\end{align}
where $ D^{\phi}_t X_t $ is the Mallivian derivative of $X_t$, $ \mathrm{d}W$ indicates the usual It\^o integral and $ \mathrm{d}^{\diamond}B^H$ denotes the FWIS integral.
\end{lemma}
		
\begin{lemma}\label{stomixed}
Let $ X_t $ be the solution of (\ref{mixSDE}). If $ f(X_t), g(X_t), h(X_t) $ be stochastic processes such that the assumptions of Lemma \ref{itomixed} are satisfied. Suppose that $ F(t,x) $ is any function satisfying continuously twice differentiable in $x$ and once in $t$. Then for $ t\in [0,T] $
\begin{align}\label{stofor}
\mathrm{d}F(t,X_t)=&\big\{\frac{\partial F}{\partial t}(t,X_t)+f(X_t)\frac{\partial F}{\partial x}(t,X_t)\big\}\mathrm{d}t\cr
&+g(X_t)\frac{\partial F}{\partial x}(t,X_t)\mathrm{d}^{\circ}W_t+h(X_t)\frac{\partial F}{\partial x}(t,X_t) \mathrm{d}^{\circ}B^H_t, \quad a.s.
\end{align}
\end{lemma}

\section{Derivation of Eq. (\ref{Malfbmx-1})}\label{app-a1}
Fixed $ r $, from (\ref{mixSDE-add}), using Lemma \ref{stomixed} and the fact that $ D^{\phi}_rX_0=0 $, one has
\begin{align}\label{Malz}
\mathrm{d}D^{\phi}_rX_t=& D^{\phi}_rf(X_t)\mathrm{d}t+D^{\phi}_rg(X_t)\mathrm{d}^{\circ}W_t+\sigma_B \phi(r,t)\mathrm{d}t\cr
=&
f'(X_t)D^{\phi}_rX_t\mathrm{d}t+g'(X_t)D^{\phi}_rX_t\mathrm{d}^{\circ}W_t+\sigma_B\phi(r,t)\mathrm{d}t.
\end{align}
			
Suppose the solution have the product form
\begin{align}\label{a1-1}
\mathrm{d}D^{\phi}_rX_t = S_tZ_t, \quad
\mathrm{d}S_t= f'(X_t)S_t\mathrm{d}t+g'(X_t)S_t\mathrm{d}^{\circ}W_t,\quad S_0=1,
\end{align}
and
\begin{align}\label{a1-2}
\mathrm{d}Z_t=A_t\mathrm{d}t+C_t\mathrm{d}^{\circ}W_t,\quad Z_0=D^{\phi}_rX_0=0,
\end{align}
where the (perhaps random) processes $ A $ and $ C $ are to be select. Then applying the It\^o formula for the Stratonovich differential, we have
\begin{align*}
\mathrm{d}D^{\phi}_rX_t&= \mathrm{d}(S_tZ_t)=S_t\mathrm{d}Z_t+Z_t\mathrm{d}S_t\cr
&=S_t\big(A_t\mathrm{d}t+C_t\mathrm{d}^{\circ}W_t\big)+Z_t\big(f'(X_t)S_t\mathrm{d}t+g'(X_t)S_t\mathrm{d}^{\circ}W_t\big)\cr
&=A_tS_t\mathrm{d}t+C_tS_t\mathrm{d}^{\circ}W_t+f'(X_t)D^{\phi}_rX_t\mathrm{d}t+g'(X_t)D^{\phi}_rX_t\mathrm{d}^{\circ}W_t,
\end{align*}
			
Now, choose $ A, C $ so that $ A_tS_t\mathrm{d}t+C_tS_t\mathrm{d}^{\circ}W_t=\sigma_B\phi(r,t)\mathrm{d}t$, 
and this identity will hold if $ C\equiv 0 $ and $ A_t:= \sigma_B\phi(r,t)S_t^{-1} $ are taken.
			
Next, the solution of (\ref{a1-1}) will be proven. The It\^o formula for $ \ln S_t $ yields 
\begin{align}\label{a1-3}
\mathrm{d}\ln S_t&=\frac{1}{S_t} \mathrm{d}S_t=\frac{1}{S_t}\big(f'(X_t)S_t\mathrm{d}t+g'(X_t)S_t\mathrm{d}^{\circ}W_t\big)\cr
&=f'(X_t)\mathrm{d}t+g'(X_t)\mathrm{d}^{\circ}W_t.
\end{align} 

Consequently, in view of the relationship between the It\^o integral and the Stratonovich integral, the following result is obtained
\begin{align*}
S_t&=\exp\big\{\int_{0}^{t}f'(X_u)\mathrm{d}u+\int_{0}^{t}g'(X_u)\mathrm{d}^{\circ}W_u\big\}\cr
&=\exp\big\{\int_{s}^{t}\big(f'(X_u)-\frac{1}{2}g(X_u)g''(X_u)\big)\mathrm{d}u+\int_{s}^{t}g'(X_u)\mathrm{d}W_u\big\}.
\end{align*}

It follows that $ S_t > 0 $ a.s. Substituting this into \ref{a1-2}, the following expression is obtained
\begin{align*}  
Z_t &= \int_0^t \sigma_B \phi(r,s) S_s^{-1} \mathrm{d}s.  
\end{align*}
			
By employing this result along with the aforementioned expression for $ S_t $, it can be concluded that  
\begin{align}\label{a1-4}
D^{\phi}_rX_t &= S_tZ_t=S_t\int_0^t\sigma_B\phi(r,s)S_s^{-1}\mathrm{d}s\cr
&=\sigma_B\int_{0}^{t}\exp\big\{\int_{s}^{t}\big(f'(X_u)-\frac{1}{2}g(X_u)g''(X_u)\big)\mathrm{d}u+\int_{s}^{t}g'(X_u)\mathrm{d}W_u\big\} \phi(r,s)\mathrm{d}s.
\end{align}

Eq. (\ref{Malfbmx-1}) is obtained by replacing $r$ by $t$.

\section{The memory kernel dependence term under approximation}
For the cases that $ g $ and $ h $ satisfy the commutativity condition, the memory kernel dependence term $ \Psi(x,t) $ can be reduced to
\begin{align}\label{eq3-0}
\Psi(x,t)=
\mathbb{E}\big[\int_{0}^{t}\exp\big\{\int_{s}^{t}\varphi_1(X_u)\mathrm{d}u\big\}\phi(t,s)\mathrm{d}s\mid X_t=x\big],
\end{align}
with $ \varphi_1(x)=f'(x)-\frac{h'(x)}{h(x)}f(x) $. For the term $\exp\big\{\int_{s}^{t}\varphi_1(X_u)\mathrm{d}u\big\} $ in Eq. (\ref{eq3-0}), using the VADA by Mamis et al. \cite{mamis2019systematic}, which provides an approximated by an appropriate stochastic Volterra-Taylor functional series expansion around a certain transient response moment.  A brief steps are as follows.

Firstly, decomposing $ \varphi_1 $ into its mean and fluctuating part:
\begin{align}\label{eq3-1}
\exp\big\{\int_{s}^{t}\varphi_1(X_u)\mathrm{d}u\big\}
&=\exp\big\{\int_{s}^{t}\big(\varphi_1(X_u)-\mathbb{E}\big[\varphi_1(X_u)\big]\big)\mathrm{d}u\big\}
\exp\big\{\int_{s}^{t}\mathbb{E}\big[\varphi_1(X_u)\big]\mathrm{d}u\big\}.
\end{align}

Secondly, applying a current-time approximation and taylor expansion in the first integral of the right-hand side of 
Eq. (\ref{eq3-1}), one has
\begin{align}\label{eq3-2}
\exp\big\{\int_{s}^{t}\varphi_1(X_u)\mathrm{d}u\big\}
&=\exp\big\{\big(\varphi_1(X_t)-\mathbb{E}\big[\varphi_1(X_t)\big]\big)(t-s)\big\}
\exp\big\{\int_{s}^{t}\mathbb{E}\big[\varphi_1(X_u)\big]\mathrm{d}u\big\}\cr
&=\Big\{1+\big(\varphi_1(X_t)-\mathbb{E}\big[\varphi_1(X_t)\big]\big)(t-s)
+
\frac{1}{2}\big(\varphi_1(X_t)-\mathbb{E}\big[\varphi_1(X_t)\big]\big)^2(t-s)^2\Big\}\cr
&\quad \times \exp\big\{\int_{s}^{t}\mathbb{E}\big[\varphi_1(X_u)\big]\mathrm{d}u\big\}.
\end{align}

Substituting Eq. (\ref{eq3-2}) into Eq. (\ref{eq3-0}), the memory kernel dependence term $ \Psi(x,t) $ becomes
\begin{align}\label{eq3-3}
\Psi(x,t)=&
\int_{0}^{t}\Big\{1+\big(\varphi_1(x)-\mathbb{E}\big[\varphi_1(X_t)\big]\big)(t-s)+\frac{1}{2}\big(\varphi_1(x)-\mathbb{E}\big[\varphi_1(X_t)\big]\big)^2(t-s)^2\Big\}\cr
&\quad \times \exp\big\{\int_{s}^{t}\mathbb{E}\big[\varphi_1(X_u)\big]\mathrm{d}u\big\}\phi(t,s)\mathrm{d}s\cr
=&
\int_{0}^{t}\exp\big\{\int_{s}^{t}\mathbb{E}\big[\varphi_1(X_u)\big]\mathrm{d}u\big\}\phi(t,s)\mathrm{d}s\cr
&+
\big(\varphi_1(x)-\mathbb{E}\big[\varphi_1(X_t)\big]\big)\int_{0}^{t}\exp\big\{\int_{s}^{t}\mathbb{E}\big[\varphi_1(X_u)\big]\mathrm{d}u\big\}(t-s)\phi(t,s)\mathrm{d}s\cr
&+\frac{1}{2}\big(\varphi_1(x)-\mathbb{E}\big[\varphi_1(X_t)\big]\big)^2\int_{0}^{t}\exp\big\{\int_{s}^{t}\mathbb{E}\big[\varphi_1(X_u)\big]\mathrm{d}u\big\}(t-s)^2\phi(t,s)\mathrm{d}s.
\end{align}

The response PDF is obtained numerically by combining Eq. (\ref{eq3-3}) with the memFPK equation. {
\begin{remark}
Specifically, through numerical experiments, the VADA method truncated at second order (as shown in Eqs. (C.3) and (C.4)) provides reliable accuracy for the memory kernel term through numerical experiments. This offers a appropriate positive diffusion coefficient ensuring the numerical stability. Truncation at a lower or higher order would cause the diffusion term to diverge or become negative, and significantly reducing the accuracy and leading the numerical instability. A more rigorous theoretical for approximation of systematic error analysis would further strengthen the framework.
\end{remark}}

\end{sloppypar}

\end{document}